\documentclass{amsart}

\usepackage{pstricks, pst-tree, pst-node} 
\usepackage[latin1]{inputenc}
\usepackage{amssymb}
\usepackage{xfrac}     
\usepackage{lmodern}
\usepackage{stmaryrd}
\usepackage{graphicx}
\usepackage{amsmath}  
\usepackage{amsthm}
\usepackage{wrapfig} 
\usepackage{color}  
\usepackage{microtype} 
\usepackage{setspace}
\usepackage{amsfonts}
\usepackage{wasysym}
\usepackage{cancel}
\usepackage{verbatim}
\usepackage{dsfont}
\usepackage{paralist}

\DeclareMathOperator{\p}{\varphi}
\DeclareMathOperator{\bp}{\mathfrak{b}_{\varphi}}

\DeclareMathOperator{\FF}{\mathbb{F}}

\DeclareMathOperator{\HH}{\mathbb{H}}

\DeclareMathOperator{\NN}{\mathbb{N}}

\DeclareMathOperator{\Ima}{Im}
\DeclareMathOperator{\lrp}{\langle\!\langle}
\DeclareMathOperator{\rrr}{\rangle\!\rangle}
\DeclareMathOperator{\rrp}{]]}
\DeclareMathOperator{\bb}{\mathfrak{b}}

\DeclareMathOperator{\w}{\wedge}
\DeclareMathOperator{\de}{d\!}

\DeclareMathOperator{\coker}{coker}
\DeclareMathOperator{\id}{id}

\newcommand{\cha}{characteristic }

\newtheoremstyle{plain2}
  {10pt}   
  {10pt}   
  {\itshape}  
  {0pt}       
  {\bfseries} 
  {}         
  {5pt plus 1pt minus 1pt} 
  {}          
  
 \newtheoremstyle{beweis}
  {10pt}   
  {10pt}   
  {\normalfont}  
  {0pt}       
  {\bfseries} 
  {:}         
  {5pt plus 1pt minus 1pt} 
  {}          

\newtheoremstyle{definition2}
  {10pt}   
  {10pt}   
  {\normalfont}  
  {0pt}       
  {\bfseries} 
  {}         
  {5pt plus 1pt minus 1pt} 
  {}          

\theoremstyle{plain2}
\newtheorem{satz}{Satz}[section]
\newtheorem{lem}[satz]{Lemma}
\newtheorem{pro}[satz]{Proposition}
\newtheorem{theo}[satz]{Theorem}
\newtheorem{coro}[satz]{Corollary}

\theoremstyle{definition2}
\newtheorem{defi}[satz]{Definition}
\newtheorem{rem}[satz]{Remark}
\newtheorem{bsp}[satz]{Example}

\theoremstyle{beweis}
\newtheorem*{prf}{Proof}

\begin{document}

\title[Witt and cohomology kernels for purely inseparable extensions]{The behavior of differential, quadratic and bilinear forms under purely inseparable field extensions}

\author{Marco Sobiech}
\address{Fakult\"at f\"ur Mathematik, Technische Universit\"at Dortmund, D-44221 Dortmund, Germany}
\email{marco.sobiech@tu-dortmund.de}
\date{\today}

\begin{abstract}
Let $F$ be a field of characteristic $p$ and let $E/F$ be a purely inseparable field extension. We study the group 
$H_p^{n+1}(F):=\coker(\wp: \Omega_F^{n} \to \Omega_F^n / \de \Omega_F^{n-1})$ of  classes of differential forms under the extension $E/F$ and give a system of 
generators of $H_p^{n+1}(E/F)$. In the case $p=2$, 
we use this to determine the kernel $W_q(E/F)$ of the restriction map 
$W_q(F) \to W_q(E)$ between the groups of nonsingular quadratic forms over $F$ and over $E$.
We also deduce the corresponding result for the bilinear Witt kernel $W(E'/F)$ of the restriction map 
$W(F) \to W(E')$, where $E'/F$ denotes a modular purely inseparable field extension. 
\end{abstract}



\maketitle

\begin{section}{Introduction}

One important aspect in the theory of quadratic and bilinear forms is to determine the behavior under field extensions, 
i.e. we want to determine which nonsingular quadratic form becomes hyperbolic and which symmetric bilinear form becomes metabolic under a given field extension. 
The main goal of this work is to classify these forms for purely inseparable extensions by using differential forms. 
For $H_p^{n+1}(F)=\coker(\wp:\Omega_F^n \to \Omega_F^n/\de\Omega_F^{n-1})$ and 
$\nu_n(F)=\ker(\wp:\Omega_F^n \to \Omega_F^n/\de\Omega_F^{n-1})$, 
we will determine generating systems for $H_p^{n+1}(E/F)$, $\Omega^n(E'/F)$ and $\nu_{n}(E'/F)$, where $E/F$ resp. $E'/F$ denotes a purely inseparable 
extension resp. a modular purely inseparable extension. These results are given in Theorem \ref{1} and Theorem \ref{41}. 
The kernels will be translated to quadratic and bilinear forms in characteristic 2 by the famous isomorphisms $H_2^{n+1}(F)\cong I^nW_q(F)/I^{n+1}W_q(F)$ and 
$\nu_n(F)\cong I^n(F)/I^{n+1}(F)$ due to Kato. With this, we will determine the quadratic resp. bilinear Witt kernel 
for the given extension, which are stated in Theorem \ref{2} and Theorem \ref{3}.

There has already been some progress in determining bilinear and quadratic Witt kernels for some given field extensions in characteristic 2. A nice short overview of some of the known Witt kernels in characteristic 2 can for example be found in \cite{Hof1}.

We note that in the case $p=2$, the same result as our main Theorem \ref{1} was computed by Laghribi, Aravire and O'Ryan using somewhat different methods independently of this work.

\end{section}

\begin{section}{A short introduction to differential forms}

We refer to \cite{Bae1}, \cite{Car} or \cite{Bae2}(for $p=2$) for any undefined terminology or any basic facts about
differential forms that we do not mention explicitly. We will assume $0\in \NN$.

Let us now start by recalling some basic facts about differential forms. Let $F$ always be a field of \cha $p>0$, if not stated otherwise. Further, let $\Omega_F^1$ be the $F$-vector space 
of absolute differential 1-forms, i.e. $\Omega_F^1$ is the $F$-vector space generated by the symbols $\de a$ with $a\in F$, together with the relations 
$\de \, (a+b)=\de a + \de b$ and $\de \, (ab)=a\de b+b\de a$ for all $a,b\in F$. Since $\de\left( F^p\right)=0$ we readily see that the map $\de: F \to \Omega_F^1 , a \mapsto \de a$ is a 
$F^p$-derivation. 

By $\Omega_F^n$ we denote the $n$-exterior power $\bigwedge^{n}\Omega_F^1$, which means $\Omega_F^n$ is a $F$-vector space generated by 
the symbols $\de a_1 \w \ldots \w \de a_n$ with $a_1,\ldots,a_n\in F$. With this, the map $\de$ can be extended to a $F^p$-linear map 
$\de: \Omega_F^{n-1} \to \Omega_F^n , a\de a_1\w \ldots \w \de a_{n-1} \mapsto \de a \w \de a_1\w \ldots \w \de a_{n-1}$. The image $\de \Omega_F^{n-1}$ of this map is called the 
set of exact forms and is additively generated by forms of the type $\de a_1 \w \ldots \w \de a_{n}$. We further set $\Omega_F^0:=F$ and $\Omega_F^n:=0$ for $n<0$.

Let us now fix a $p$-basis $\mathcal{B}=\{b_i \mid i\in I \}$ of $F$ over $F^p$, i.e. we have $F^p(\mathcal{B})=F$ and for every finite subset 
$\{b_{i_1},\ldots ,b_{i_k}\}\subset \mathcal{B}$, we have $[F^p(b_{i_1},\ldots ,b_{i_k}) : F^p]=p^k$, where the second condition is called the $p$-independence of $\mathcal{B}$. For further details on $p$-bases see \cite{Car}. We may assume that the index set $I$ is well-ordered and transfer this ordering to $\mathcal{B}$ by setting $b_i<b_j$ iff $i<j$. 
We further define the set $\Sigma_n :=\{\sigma : \{1,\ldots,n\} \to I \mid \sigma(i)<\sigma(j) \text{ for } i<j \}$ and equip it with the lexicographic ordering, i.e. 
we set $\sigma \leq \delta$ iff $\sigma(j) \leq \delta(j)$ and $\sigma(i)=\delta(i)$ for some $j\in \{1,\ldots,n\}$ and all $i < j$. Then it is well known 
that the logarithmic differential forms $\{\frac{\de b_{\sigma}}{b_{\sigma} } \mid  \sigma \in \Sigma_n \}$ form a $F$-basis \label{basis}of $\Omega_F^n$, where we set 
$\frac{\de b_{\sigma}}{b_{\sigma}}:= \frac{\de b_{\sigma(1)}}{b_{\sigma(1)}}\w\ldots \w \frac{\de b_{\sigma(n)}}{b_{\sigma(n)}}$. We further 
define $F_j$ resp. $F_{<j}$ to be the subfield of $F$ generated over $F^p$ by the elements $b_i$ with $i\leq j$ resp. $i<j$. 
We also set $\Omega_{F,\alpha}^n$ resp. $\Omega_{F,<\alpha}^n$ to be the subspace of $\Omega_F^n$ generated by the elements $\frac{\de b_{\sigma}}{b_{\sigma}}$ with 
$\sigma \leq \alpha$ resp. $\sigma < \alpha$. This gives rise to a filtration of the space $\Omega_F^n$.

After fixing a $p$-basis of $F$, we can define the additive homomorphism 
$$s_p: \Omega_F^n \to \Omega_F^n \ ,\  \sum_{\sigma}a_{\sigma}\frac{\de b_{\sigma}}{b_{\sigma}} \mapsto 
\left( \sum_{\sigma}a_{\sigma}\frac{\de b_{\sigma}}{b_{\sigma}} \right)^{[p]}:=\sum_{\sigma}a_{\sigma}^p\frac{\de b_{\sigma}}{b_{\sigma}}.$$ 
Using this, we can extend the usual Artin-Schreier map $\wp: F \to F$ to $\wp:= s_p- \id : \Omega_F^n \to \Omega_F^n$. Both maps $s_p$ and $\wp$ are obviously basis dependent, but a change of basis changes the image 
of $s_p$ and $\wp$ by an exact form, which means that 
\begin{align*}
s_p: \Omega_F^n \to \Omega_F^n / \de \Omega_F^{n-1} \ \text{ and }\  \wp: \Omega_F^n \to \Omega_F^n / \de \Omega_F^{n-1}
\end{align*} 
are well defined. With this we set 
\begin{align*}
H_p^{n+1}(F)&:=\Omega_F^n/\left( \wp \Omega_F^n + \de \Omega_F^{n-1} \right)=\coker\wp \\
\nu_n(F)&:=\ker \wp
\end{align*}
A differential form $\omega$ viewed in $H_p^{n+1}(F)$ will be denoted by $\overline{\omega}$. Note that for $n=0$, we have $H_p^1(F)=F/\wp(F)$. 

For an arbitrary field extension $L/F$, 
we define 
\begin{align*}
H_p^{n+1}(L/F)&:= \ker (  H_p^{n+1}(F) \to H_p^{n+1}(L) ) \\
\Omega^n (L/F)&:=\ker ( \Omega_F^n \to \Omega_L^n ) \\
\nu_n (L/F)&:= \ker ( \nu_n (F) \to \nu_n (L)  )= \Omega^n(L/F) \cap \nu _n(F) 
\end{align*}
 

We will often freely use the following well known result and will prove it for the reader's convenience.
\begin{lem}\label{30}
In $H_p^{n+1}(F)$ we have $\overline v = \overline{ v^{[p^i]}}$ for all $v \in \Omega_F^n $ and all $ i\in \NN$.
\end{lem}

\begin{prf}
For a $v\in\Omega_F^n$, the following equation holds
$$ v^{[p^i]}\equiv v^{[p^i]} +\wp (-v^{[p^{i-1}]}) \equiv v^{[p^{i}]} +(-1)^{p} v^{[p^{i}]}+v^{[p^{i-1}]}   \equiv  v^{[p^{i-1}]}   \! \mod (\wp\Omega_F^n + \de \Omega_F^{n-1}).  $$
An induction on $i\in\NN$ finishes the proof. \hfill $\square$
\end{prf}

Since we will start by studying the behavior of differential forms under purely inseparable field extensions, we will close this section with a well known fact for these kinds 
of field extensions. For this, let us recall that a field extension $L/F$ is called purely inseparable if, for all $\ell \in L$, there exists a $k\geq 0$ with $\ell^{p^k}\in F$. 
Using this definition together with $H_p^1(F)=F/\wp(F)$, the following lemma can be viewed as a corollary of Lemma \ref{30}.

\begin{lem}\label{355}
Let $L/F$ be a purely inseparable field extension and $x\in F$. If $x\in\wp(L)$, then $x\in\wp(F)$.
In particular, we have $H_p^1(L/F)=\{\overline{0}\}$.
\end{lem}

\end{section}

\begin{section}{The kernel $H_p^{n+1}(E/F)$ for purely inseparable field extensions}\label{sec3}

In this section, our main goal is to compute the kernel $H_p^{n+1}(E/F)$ where $E$ is a finite purely inseparable field extension of $F$. For that let us recall that each purely 
inseparable field extension of $F$ can be realized as $E=F(\sqrt[p^{m_1}]{b_1}, \ldots , \sqrt[p^{m_r}]{b_r})$ with integers $m_1,\ldots,m_r \in \NN_{>0}$ and $b_1,\ldots,b_r \in F$. 
The exponent the the extension $E/F$ is defined as $e:=\exp(E/F):=\min \{k \in \NN \mid x^{p^k}\in F$ for all $x\in E \}$. Note that for $b_1,\ldots,b_r \in F\setminus F^p$, we have 
$e=\max\{m_1,\ldots,m_r\}$. To find generators of $H_p^{n+1}(E/F)$, we will need the following two cases as intermediate steps
\begin{itemize}
\item $r=1$ and $\exp(E/F)=1$,
\item $r$ arbitrary and $\exp(E/F)=1$,
\end{itemize}
which will be discussed in \ref{3.3} and \ref{34}. These results are known for $p=2$ and were proven in \cite{Bae2} and \cite{Lag2}. 
For the generalization of the known results to an arbitrary prime $p$, 
we will use the following lemmas.

\begin{lem}\label{31}{{\rm(\cite{Bae1}, Lemma 3.1)}}
Let $\mathcal{B}=\{b_i \mid i\in I \}$ be a $p$-basis of the field $F$. Assume 
$\sum_{\sigma \leq \alpha }c_{\sigma} \frac{\de b_{\sigma}}{b_{\sigma}}\in \de \Omega_F^{n-1},$
with $c_{\alpha}\neq 0$. Then there exist elements $M_{ij}\in F_{<\alpha(i)}$ with $1 \leq i\leq n \ , \ 1 \leq j \leq p-1$, so that 
$$c_{\alpha} = \sum_{i=1}^n\sum_{j=1}^{p-1}M_{ij}b_{\alpha(i)}^j.$$
\end{lem}

The next lemma is a weaker version of the famous lemma of Kato. The original Kato's lemma \cite[Lemma 2]{Kat1} can only be used for fields which are $p-1$ closed, i.e. $F^{p-1}=F$. The following weaker version does not need this assumption and will suffice for all the upcoming computations in this section.

\begin{lem}\label{32}{{\rm(\cite{Bae1}, Lemma 3.2)}}
Let $\mathcal{B}=\{b_i \mid i\in I \}$ be a $p$-basis of the field $F$, $\alpha\in \Sigma_n$ and $u\in F$ such that 
$$\wp(u)\frac{\de b_{\alpha}}{b_{\alpha}}\in \de \Omega_F^{n-1} + \Omega_{F,<\alpha}^n.$$
Then we have 
$$ u\frac{\de b_{\alpha}}{b_{\alpha}} = \sum_{(i_1,\ldots,i_n)} \frac{\de a_{i_1}}{a_{i_1}} \w \ldots \w \frac{\de a_{i_n}}{a_{i_n}} + v,$$
with suitable $a_{i_j}\in F$ for $1 \leq j \leq n$ and $v\in \Omega_{F,<\alpha}^n$.
\end{lem}

With Lemma \ref{32}, it is easy to show that $\nu_n(F)$ is additively generated by the logarithmic differential forms $\frac{\de a_1}{a_1}\w \ldots \w \frac{\de a_n}{a_n}$ with 
$a_1,\ldots,a_n\in F$. Since we will need 
the original version of Kato's Lemma in the last section, we will state it here for easy reference.

\begin{lem}\label{02}{{\rm(\cite{Kat1}, Lemma 2)}}
Let $\mathcal{B}=\{b_i \mid i\in I \}$ be a $p$-basis of the field $F$ with $F^{p-1}=F$, $\alpha\in \Sigma_n$ and $u\in F$ such that 
$$\wp(u)\frac{\de b_{\alpha}}{b_{\alpha}}\in \de \Omega_F^{n-1} + \Omega_{F,<\alpha}^n.$$
Then we have 
$$ u\frac{\de b_{\alpha}}{b_{\alpha}} =  \frac{\de a_{1}}{a_{1}} \w \ldots \w \frac{\de a_{n}}{a_{n}} + v,$$
with suitable $a_i \in F_{\alpha(i)}^*$ for $1 \leq i \leq n$ and $v\in \Omega_{F,<\alpha}^n$.
\end{lem}

Note that the condition $F^{p-1}=F$ is obviously true for all fields of characteristic $2$.

\begin{rem}\label{00}
In most of the upcoming proofs and computations, we will have one single element $\overline{\omega}\in H_p^{n+1}(L/F)$ for some finitely generated algebraic field extension 
$L=F(\beta_1,\ldots,\beta_r)/F$. Since the minimal polynomial over $F$ of each $\beta_i$, the form $\omega$ and the relation 
$\omega_L= \wp(u) + \de v$ for some $u\in\Omega_L^n$ and $v\in\Omega_L^{n-1}$ only need a finite number of coefficients of $F$ to be described, we may assume that $\omega$ is defined 
over a field $F'$, which is generated by finitely many elements over $\FF_p$. By setting $L'=F'(\beta_1,\ldots,\beta_r)$ we get that $\overline{\omega}_{L'}=0$, which means 
$\overline{\omega}\in H_p^{n+1}(L'/F')$. Following these ideas, we may always assume 
that $\omega$ is defined over a field with a finite $p$-basis. Note that 
the polynomials of the $\beta_i$ over $F'$ are the same as over the field $F$.
\end{rem}

\begin{rem}\label{04}
(i) Let $\mathcal{B}=\{c_i \mid i\in I\}$ be a $p$-basis of $F$ with well-ordered indexset $I$. For a fixed $s\in \NN$, let $i_1,\ldots,i_s\in I$ be the smallest elements in $I$ 
with $i_1 < \ldots < i_s$. By this we have $\Sigma_n=\Sigma_{n,>s}\mathbin{\dot{\cup}}\Sigma_{n\leq s}$ with 
\begin{align*}
\Sigma_{n,>s}&:= \{ \pi\in \Sigma_n \mid  \Ima(\pi)\cap \{i_1,\ldots,i_s\} = \varnothing  \}  \\
\Sigma_{n,\leq s}&:= \{ \pi\in \Sigma_n \mid  \Ima(\pi)\cap \{i_1,\ldots,i_s\} \neq \varnothing  \}  
\end{align*} 
Thus, for every $\omega \in \Omega_F^n$, we can find unique forms $\omega_{>s},\omega_{\leq s}\in \Omega_F^n$ with $\omega=\omega_{> s} +\omega_{\leq s}$ where in the basis expansion 
of $\omega_{>s}$ none of the summands is divided by any of $\de c_{i_1}, \ldots, \de c_{i_s}$ and in the basis expansion of $\omega_{\leq s}$ every summand is divided by at least one of the forms 
$\de c_{i_1}, \ldots, \de c_{i_s}$. If there is no confusion about the number $n$, we omit it as an index.

(ii) Let $b_1,\ldots,b_r\in F$ be $p$-independent and $\mathcal{B}$ a $p$-basis of $F$ containing $b_1,\ldots,b_r$. 
An easy computation shows that a $p$-basis of the purely inseparable field extension 
$E=F\left(\sqrt[p^{m_1}]{b_1},\ldots,\sqrt[p^{m_r}]{b_r} \right)$, is given by $\mathcal{B}_E=(\mathcal{B}\setminus\{b_1,\ldots,b_r\})\cup \{\sqrt[p^{m_1}]{b_1},\ldots,\sqrt[p^{m_r}]{b_r}\}$. 
Thus if we set $U:=\{\sqrt[p^{m_1}]{b_1},\ldots,\sqrt[p^{m_r}]{b_r}\}$, we get 
$$\Omega_E^n= E \left(\Omega_F^n\right)_E \oplus   \bigoplus_{\substack{\varnothing \neq \{a_1,\ldots,a_{\ell}\}\subseteq U\\ 1\leq \ell \leq r}}
E \de  a_1  \w \ldots \w \de a_{\ell}  \w \left(\Omega_F^{n-\ell}\right)_E.$$
If the $\sqrt[p^{m_1}]{b_1},\ldots,\sqrt[p^{m_r}]{b_r}$ are the first $r$ elements of $\mathcal{B}_E$, this decomposition coincides with the one given above.
\end{rem}

In \cite[Lemma 7.1]{Hof1} Hoffmann and Dolphin showed that for the extension $F(\sqrt[p^m]{b})/F$ with $b\in F\setminus F^p$, we have 
\begin{equation}\label{eq0}
\Omega^n(F(\sqrt[p^m]{b})/F)=\de b \w \Omega_F^{n-1}.
\end{equation}
This will be used to prove the following

\begin{pro} \label{3.3}
For $b\in F$ and $E:=F\left( \sqrt[p]{b} \right)$, we have 
$$ H_p^{n+1}(E/F)= \overline{  \de b \w \Omega_F^{n-1}  }= \overline{  b\de\left( \Omega_F^{n-1}\right)}  . $$
\end{pro}

\begin{prf}
For $b\in F^p$, we have $E=F$, thus $H_p^{n+1}(E/F)=\{0\}=\overline{\de b \w\Omega_F^{n-1} }$. So we can assume $b\in F\setminus F^p$. 

First we have to check that each form of the type $\overline{ \de b \w y }$ with $y\in \Omega_F^{n-1}$ is an element of $H_p^{n+1}(E/F)$. But this is obviously true, since we have 
$b\in E^p$. 

Let us now prove the reverse inclusion. Set  $\beta:=\sqrt[p]{b}$. Since $b\in F\setminus F^p$, we may choose $b$ to be the first element of a $p$-basis $\mathcal{B}=\{b_i \mid i\in I\}$ of $F$. With this, we get a $p$-basis of $E$ by 
$\mathcal{B}_E:=(\mathcal{B}\setminus \{b\})\cup \{\beta\}$. Without loss of generality, we can choose $\beta$ to be the last element of $\mathcal{B}_E$,
i.e. $b'<\beta$ for all $b'\in \mathcal{B}\setminus \{b\}$.

Now choose an element $\overline{\omega}\in H_p^{n+1}(E/F)$. Since $\overline{\omega_E} = 0 \in H_p^{n+1}(E)$, we can find $u\in \Omega_E^n$ and $v \in \Omega_E^{n-1}$ with 
$\omega_E = \wp(u) + \de v.$ Using the basis described above, let $\gamma$ be the maximal multiindex in the basis expansion of $\omega$, i.e. 
$\omega = \sum_{\sigma \leq \gamma} \omega_{\sigma} \frac{\de b_{\sigma}}{b_{\sigma}}$. Since we already know that $\de b =0$ in $E$, by Remark \ref{04} 
we may assume that the basis expansion of $\omega$ is free of forms containing $\de b$, which means we have $\sigma \in \Sigma_{n,>1}$ for all the 
multiindices $\sigma$ used to describe $\omega$. With this 
assumption, this representation of $\omega$ is also the unique representation of $\omega_E$ using the $p$-basis $\mathcal{B}_E$ over the field $E$. 
Let $\delta$ be the maximal 
multiindex in the basis expansion of $u$, i.e. 
$u = \sum_{\sigma \leq \delta} u_{\sigma} \frac{\de b_{\sigma}}{b_{\sigma}}$. Using Remark \ref{00}, we may assume that there are only finitely many multiindices lower than $\gamma$ resp. 
$\delta$.

Assume $\delta > \gamma$. Then we have 
\begin{align*}
0 \equiv \wp(u_{\delta})\frac{\de b_{\delta}}{b_{\delta}} + \de v \mod \Omega_{E,<\delta}^n.
\end{align*}
We can now apply Lemma \ref{32} to get
\begin{align*}
u_{\delta}\frac{\de b_{\delta}}{b_{\delta}} = \sum_{(i_1,\ldots,i_n)} \frac{\de a_{i_1}}{a_{i_1}} \w \ldots \w \frac{\de a_{i_n}}{a_{i_n}} + u_{<\delta}
\end{align*}
with some suitable $a_{i_j}\in E$ and $u_{<\delta} \in \Omega_{E,<\delta}^n$. Since $u$ can be written as $u=u_{\delta}\frac{\de b_{\delta}}{b_{\delta}} + u_{<\delta}'$ with
$u_{<\delta}'\in \Omega_{E,<\delta}^n$, we 
get that $u= \sum_{(i_1,\ldots,i_n)} \frac{\de a_{i_1}}{a_{i_1}} \w \ldots \w \frac{\de a_{i_n}}{a_{i_n}} + u_{<\delta}' +u_{<\delta}$. Since all the forms 
 $\frac{\de a_{i_1}}{a_{i_1}} \w \ldots \w \frac{\de a_{i_n}}{a_{i_n}}$ are elements of $\nu_n(E)$, we see that 
$\wp(u)=\wp(u_{<\delta}'')$, with $u_{<\delta}'':=u_{<\delta}' +u_{<\delta} \in \Omega_{E,<\delta}^n$. Thus we can lower the maximal multiindex of $u$ and because of Remark \ref{00}, 
after a finite number of repetitions, we get $\delta \leq \gamma$. 

So let us now assume $\delta \leq \gamma$. With this we have 
\begin{align*}
\omega_{\gamma}\frac{\de b_{\gamma}}{b_{\gamma}} \equiv \wp(u_{\gamma})\frac{\de b_{\gamma}}{b_{\gamma}} + \de v \mod \Omega_{E,<\gamma}^n,
\end{align*}
which means 
\begin{align*}
(\omega_{\gamma}- \wp(u_{\gamma})) \frac{\de b_{\gamma}}{b_{\gamma}} \equiv  \de v  \mod \Omega_{E,<\gamma}^n.
\end{align*}
Applying Lemma \ref{31}, we find $M_{ij}\in E_{<\gamma(i)}=E^p\left(b'\in\mathcal{B}_E  \mid b' < b_{\gamma(i)} \right)$ with 
\begin{align}\label{eq1}
(\omega_{\gamma}- \wp(u_{\gamma}))= \sum_{i=1}^n\sum_{j=1}^{p-1} M_{ij} b_{\gamma(i)}^j.
\end{align}
Since by assumption both $\omega$ and $\omega_E$ have the same basis representation, we have $b_{\gamma(1)}, \ldots , b_{\gamma(n)}\in F$ and since $E^p\subset F$, 
we also have $M_{ij}\in F$ for all suitable $i,j$, hence 
$\sum_{i=1}^n\sum_{j=1}^{p-1} M_{ij} b_{\gamma(i)}^j\in F$. Since $u_{\gamma}\in E$, we find $u_0,\ldots,u_{p-1}\in F$ with $u_{\gamma}=\sum_{\ell=0}^{p-1}\beta^{\ell} u_{\ell}$. 
By Equation \eqref{eq1}, we have $\wp(u_{\gamma})\in F$, so by using Lemma \ref{355}, we may assume $u_{\gamma}=u_0\in F$.
With this, we get
\begin{align*}
(&\omega_{\gamma}-\wp(u_0))\frac{\de b_{\gamma}}{b_{\gamma}} =
  \sum_{i=1}^n\sum_{j=1}^{p-1} M_{ij} b_{\gamma(i)}^j \frac{\de b_{\gamma(1)}}{b_{\gamma(1)}} \w \ldots \w \frac{\de b_{\gamma(n)}}{b_{\gamma(n)}} \\
  &= \sum_{i=1}^n\sum_{j=1}^{p-1} \Biggl(  \de \left(  (-1)^{i+1}\frac{M_{ij}}{j} b_{\gamma(i)}^j
  \frac{\de b_{\gamma(1)}}{b_{\gamma(1)}} \w \ldots \w\frac{\de b_{\gamma(i-1)}}{b_{\gamma(i-1)}} \w \frac{\de b_{\gamma(i+1)}}{b_{\gamma(i+1)}}\w \ldots \w \frac{\de b_{\gamma(n)}}{b_{\gamma(n)}}
  \right)\\ 
  &\qquad \quad -  \frac{M_{ij}}{j} b_{\gamma(i)}^j  
  \frac{\de b_{\gamma(1)}}{b_{\gamma(1)}} \w \ldots \w \frac{\de b_{\gamma(i-1)}}{b_{\gamma(i-1)}} \w \frac{\de M_{ij}}{M_{ij}}\w 
  \frac{\de b_{\gamma(i+1)}}{b_{\gamma(i+1)}}\w \ldots \w\frac{\de b_{\gamma(n)}}{b_{\gamma(n)}} \Biggl) \\
  &\equiv \de t \mod \Omega_{F,<\gamma}^n
\end{align*}
where
$$t:=\sum_{i=1}^n\sum_{j=1}^{p-1}  \left( (-1)^{i+1} \frac{M_{ij}}{j} b_{\gamma(i)}^j
  \frac{\de b_{\gamma(1)}}{b_{\gamma(1)}} \w \ldots \w\frac{\de b_{\gamma(i-1)}}{b_{\gamma(i-1)}} \w \frac{\de b_{\gamma(i+1)}}{b_{\gamma(i+1)}}\w \ldots \w\frac{\de b_{\gamma(n)}}{b_{\gamma(n)}}
  \right),$$ 
 with $t\in \Omega_F^{n-1}$ and because $M_{ij}\in E_{<\gamma(i)}$. So now we have 
\begin{align}\label{eq2}
\omega':=(\omega_{\gamma}-\wp(u_0))\frac{\de b_{\gamma}}{b_{\gamma}} - \de t \in \Omega_{F,<\gamma}^n
\end{align}
By inserting Equation \eqref{eq2} in the decomposition $\omega= \omega_{\gamma}\frac{\de b_{\gamma}}{b_{\gamma}} + \omega_{<\gamma}$ with $\omega_{<\gamma}\in\Omega_{F,<\gamma}^n$, we get 
\begin{align*}
\omega \ =\  \omega_{\gamma}\frac{\de b_{\gamma}}{b_{\gamma}} + \omega_{<\gamma} 
\ =\ \wp(u_0) \frac{\de b_{\gamma}}{b_{\gamma}} + \de t + \omega' + \omega_{<\gamma} 
\ =\  \wp(u_0) \frac{\de b_{\gamma}}{b_{\gamma}} + \de t + \omega''
\end{align*}
with $\omega''=\omega' + \omega_{<\gamma}\in\Omega_{F,<\gamma}^n $. This shows $\overline{\omega}=\overline{\omega''}$, thus we can also lower the maximal multiindex of the 
form $\omega$ over $E$ without changing its class in $H_p^{n+1}(E)$. Repeating this together with lowering 
the maximal multiindex of $u$ whenever possible and by using \ref{00}, 
after a finite number of reductions we may assume $\omega_E=0 \in \Omega_E^n$. Under the assumption we made on $\omega$ and by using Equation \eqref{eq0}, we 
obtain $\omega=0$. With this we readily get $H_p^{n+1}(E/F)\subset \overline{\de b \w \Omega_F^{n-1} }$.

The second description of the kernel can be easily seen by $\overline{\de b \w u}= \overline{ b \de\left( -u \right)}$ for all $b\in F$ and $u\in\Omega_F^{n-1}$.
 \hfill$\square$
\end{prf}

Since we now have found the kernel for purely inseparable field extensions of degree $p$, our next step is to generalize this result to field extensions of  type 
$F\left( \sqrt[p]{b_1},\ldots, \sqrt[p]{b_r} \right)/F$ with $b_1,\ldots,b_r\in F$.

\begin{theo}\label{34}
For $b_1,\ldots,b_r\in F$ and $E:=F\left( \sqrt[p]{b_1},\ldots,\sqrt[p]{b_r} \right)$, we have 
$$ H_p^{n+1}(E/F)= \sum_{i=1}^r \overline{  \de b_i \w \Omega_F^{n-1}  }=\sum_{i=1}^r \overline{  b_i  \de \Omega_F^{n-1}  } . $$
\end{theo}

\begin{prf}

Since $b_1,\ldots,b_r \subset E^p$ we obviously get $\sum_{i=1}^r \overline{  \de b_i \w \Omega_F^{n-1}  } \subset H_p^{n+1}(E/F)$. To check the reverse inclusion,  
we use induction on $r$. The case $r=1$ is covered by Proposition \ref{3.3}, so let us assume $r>1$. First assume $[E:F]=p^r$, i.e. the elements $b_1,\ldots,b_r$ are $p$-independent and set $M:=F(\sqrt[p]{b_1})$ . Hence we can choose a $p$-basis $\mathcal{B}_F=\{c_i \mid i\in I\}$ of $F$ 
with $c_{i_1}=b_1,\ldots,c_{i_r}=b_r$ where $i_1,\ldots,i_r$ denote the first $r$ elements in $I$ with $i_1< \ldots < i_r$. By Remark \ref{00}, we may assume $|I|<\infty$. 
Similar to the proof of \ref{3.3}, we then obtain a $p$-basis of $M$ by  $\mathcal{B}_M:=(\mathcal{B}_F \setminus \{b_1\}) \cup \{ \sqrt[p]{b_1}\}$. Further 
we set $\sqrt[p]{b_1}>c_i$ for all $i\in I\setminus \{i_1\}$.

Now choose an element $\overline{ \omega } \in H_p^{n+1}(E/F)$. Then we have $\overline{ \omega_M } \in H_p^{n+1}(E/M)$ and since $E=M\left(\sqrt[p]{b_2}, \ldots, \sqrt[p]{b_r} \right)$, 
we get by the induction hypothesis 
\begin{align}\label{eq3}
\omega_M = \wp(u) + \de v + \sum_{i=2}^r \de b_i \w x_i
\end{align}
for suitable $u\in \Omega_M^n$ and $v,x_i \in \Omega_M^{n-1}$. Let $\delta$ be the maximal multiindex in Equation \eqref{eq3} using the basis $\mathcal{B}_M$ 
and assume $\delta(1)>i_r$. Then we have 
\begin{align*}
\omega_{\delta} \frac{\de c_{\delta}}{c_{\delta}} \equiv \wp(u_{\delta})\frac{\de c_{\delta}}{c_{\delta}} + \de v  
\equiv \wp(u_{\delta})\frac{\de c_{\delta}}{c_{\delta}} + t \frac{\de c_{\delta}}{c_{\delta}}   \mod \Omega_{M,<\delta}^n.
\end{align*}
with $\de v - t \frac{\de c_{\delta}}{c_{\delta}} = v' \in \Omega_{M,<\delta}^{n}$. 

First suppose $c_{\delta(n)}=\sqrt[p]{	b_1}:=\beta$. Then obviously $\omega_{\delta}=0$, since $\omega$ is defined over $F$. Hence we have 
\begin{align*}
\wp(u_{\delta})\frac{\de c_{\delta}}{c_{\delta}} \equiv   \de\left(- v\right) + v' \in \de\Omega_M^{n-1} + \Omega_{M,<\delta}^n.
\end{align*}
Using Lemma \ref{32}, we get 
$$u_{\delta}\frac{\de c_{\delta}}{c_{\delta}} = \sum_{(j_1,\ldots,j_n)} \frac{\de a_{j_1}}{a_{j_1}} \w \ldots \w \frac{\de a_{j_n}}{a_{j_n}} + u',$$
with $a_{j_1},\ldots,a_{j_n}\in M$ and $u'\in\Omega_{M,<\delta}^{n}$. By writing $u= u_{\delta}\frac{\de c_{\delta}}{c_{\delta}} +u''$ with $u''\in\Omega_{M,<\delta}^{n}$, we get 
$\wp(u)=\wp(u'+u'')$, thus we can lower the maximal multiindex of $u$. Iterating this process whenever possible, we may assume that the Equation \eqref{eq3} does not contain a slot 
$\de \beta$, i.e. $c_{\delta(n)} < \beta$. 

As a next step, we will show that we may assume $v\in \Omega_F^{n-1}$. Say 
$$v = \sum_{i=0}^{p-1} \beta^i v_i + \sum_{i=0}^{p-1} \beta^i  v'_i \w \de \beta\ \ , \ \ \text{ with } v_i\in\Omega_F^{n-1} \text{ and } v'_i\in\Omega_F^{n-2}.$$
With this we get 
$$ \de v =  \sum_{i=0}^{p-1} \beta^i \de v_i +  \sum_{i=1}^{p-1} i(-1)^{n-1} \beta^{i-1}   v_i \w \de \beta +  \sum_{i=0}^{p-1} \beta^i \de v'_i \w \de \beta.$$
Since Equation \eqref{eq3} does not contain a slot $\de \beta$ and $\beta$ is part of the $p$-basis of $M$, we get 
\begin{align*}
0&=\sum_{i=1}^{p-1} i(-1)^{n-1} \beta^{i-1}   v_i   +  \sum_{i=0}^{p-1} \beta^i \de v'_i \\
 &= \sum_{i=1}^{p-1} \Big( \beta^{i-1}\left( (-1)^{n-1}i v_i + \de v'_{i-1} \right) \Big) + \beta^{p-1}\de v'_{p-1},
\end{align*}
hence  $\de v'_{p-1}=0$ and $v_i= \de\left(\frac{(-1)^n}{i}  v'_{i-1} \right)\in \de \Omega_F^{n-2}$ for $i=1,\ldots, p-1$. With this, we have 
$$\de v = \de v_0 + \sum_{i=1}^{p-1} \beta^i \de v_i = \de v_0 + \sum_{i=1}^{p-1} \beta^i \de \left(\de\left(  (-1)^ni^{-1}v'_{i-1}  \right)\right) = \de v_0.$$
So we may assume $v\in \Omega_F^{n-1}$ and $t \frac{\de c_{\delta}}{c_{\delta}}, v' \in\Omega_F^n$ as well. 

Now with $\omega_{\delta},t\in F$, we get $\wp(u_{\delta})\in F$, thus $u_{\delta}\in F$ by Lemma \ref{355}. Substituting $\omega_{\delta} \frac{\de c_{\delta}}{c_{\delta}}$ 
in $\omega=\omega_{\delta} \frac{\de c_{\delta}}{c_{\delta}} + \omega' $ with $\omega'\in\Omega_{F,<\delta}^{n-1}$, we get 
\begin{align*}
\omega &= \omega_{\delta} \frac{\de c_{\delta}}{c_{\delta}} + \omega' \\
&= \wp(u_{\delta}) \frac{\de c_{\delta}}{c_{\delta}} + t \frac{\de c_{\delta}}{c_{\delta}} +\omega' \\
&= \wp(u_{\delta}) \frac{\de c_{\delta}}{c_{\delta}} + \de v - v' +\omega' \\
&\equiv \omega'-v' \mod (\wp\Omega_F^n + \de \Omega_{F}^{n-1})
\end{align*}
with $\omega'-v' \in \Omega_{F,<\delta}^{n-1}$. Iterating this process, after a finite number of steps, 
we find  $\omega''\in \Omega_F^n$ with $\overline{\omega}=\overline{\omega''}\in H_p^{n+1}(F)$, 
so that the 
maximal multiindex of $\omega''$ is given by $\varepsilon$ with $\varepsilon (1) \leq i_r$. This means $\omega'' \in \sum_{i=2}^r \de b_i \w \Omega_F^{n-1}$, 
thus we can find forms $y_2,\ldots,y_r \in \Omega_F^{n-1}$, 
such that over the field $M$ we have
\begin{align*}
\overline{  \omega } - \sum_{i=2}^r  \overline{ \de b_i \w y_i  } =0 \in H_p^{n+1}(M).
\end{align*}
Using \ref{3.3} and by choosing a suitable $y_1\in \Omega_F^{n-1}$, we get
\begin{align*}
\overline{  \omega }- \sum_{i=2}^r \overline{ \de b_i \w y_i  }=  \overline{ \de b_1 \w y_1  },
\end{align*}
from which we finally obtain
$$\overline{ \omega  } =   \sum_{i=1}^r \overline{ \de b_i \w y_i  } \in H_p^{n+1}(F). $$
Note that we have $\omega=\omega''$, if $\delta(1)\leq i_r$ in Equation \eqref{eq3}.

If the elements $b_1,\ldots,b_r$ are $p$-dependent, i.e. $[F\left( \sqrt[p]{b_1},\ldots,\sqrt[p]{b_r} \right):F]=p^s < p^r$, we can find $\{\ell_1,\ldots,\ell_s\}\subset \{1,\ldots, r$\} with 
$F\left( \sqrt[p]{b_1},\ldots,\sqrt[p]{b_r} \right)=F\left( \sqrt[p]{b_{\ell_1}},\ldots,\sqrt[p]{b_{\ell_s}} \right)$ and the elements $b_{\ell_1},\ldots,b_{\ell_s}$ are $p$-independent. 
Applying the first part of the proof, we get
$$H_p^{n+1}(F\left( \sqrt[p]{b_1},\ldots,\sqrt[p]{b_r} \right)= \sum_{k=1}^s \overline{\de b_{\ell_k}\w \Omega_F^{n-1}} \overset{(1)}{=}\sum_{i=1}^r \overline{  \de b_i \w \Omega_F^{n-1}  }$$
where step (1) is possible because of $\de b_1,\ldots,\de b_r \in \langle \de b_{\ell_1},\ldots , \de b_{\ell_s} \rangle$ due to $b_1\ldots,b_r \in F^p(b_{\ell_1},\ldots,b_{\ell_s})$. 

\hfill$\square$
\end{prf}

To generalize \ref{34} to a purely inseparable extension of arbitrary exponent, we will need a few more results and notations.

\begin{lem}\label{366}
For $v\in \Omega_F^{n-1}$, $b_1,\ldots,b_r\in F$ and $k_1,\ldots,k_r\geq 1$ the following equation holds in $H_p^{n+1}(F)$
$$ \overline{  b_1^{k_1}\ldots b_r^{k_r} \de v  }\equiv   \sum_{j=1}^r \overline{b_j \de\left( k_j b_1^{k_1}\ldots b_{j-1}^{k_{j-1}} b_j^{k_j-1} b_{j+1}^{k_{j+1}} \ldots b_r^{k_r} v \right)}. $$

\end{lem}

\begin{prf}
For $r=1$, we will prove the lemma by induction on $k_1$. For $k_1=1$ there is nothing to prove, so assume $k_1\geq 2$ and set $b_1=b$ and $k_1=k$. Then we have
\begin{align*}
0 &\equiv \de\left(b^kv \right)\equiv b^k\de v + kb^{k-1}\de b \w  v \ \mod\left( \wp\Omega_F^n + \de \Omega_F^{n-1}\right),
\end{align*}
and we also have 
\begin{align*}
0 &\equiv k\de\left(b\left(b^{k-1}v\right) \right)\equiv kb\de\left( b^{k-1}v\right) + k\de b \w (b^{k-1}v) \ \mod\left( \wp\Omega_F^n + \de \Omega_F^{n-1}\right)
\end{align*}
from which we easily get $b^k\de v \equiv b\de\left( k b^{k-1}v\right) \mod\left( \wp\Omega_F^n + \de \Omega_F^{n-1}\right)$, since we have $k\in \FF_p \subset F^p$.

Now assume $r>1$. With this, we have 
\begin{align*}
 b_1^{k_1}\ldots b_r^{k_r} \de v &\equiv -\de\left(  b_1^{k_1}\ldots b_r^{k_r} \right) \w v \\
 &\equiv -b_1^{k_1} \de\left(b_2^{k_2}\ldots b_r^{k_r} \right) \w v - b_2^{k_2}\ldots b_r^{k_r} \de \left( b_1^{k_1} \right) \w v \\
  &\equiv - \de\left(b_2^{k_2}\ldots b_r^{k_r} \right)  \w \left(b_1^{k_1}v\right) -  \de \left( b_1^{k_1} \right) \w \left( b_2^{k_2}\ldots b_r^{k_r}v \right)\\
  &\equiv b_2^{k_2}\ldots b_r^{k_r} \de \left( b_1^{k_1} v\right)  + b_1^{k_1} \de \left( b_2^{k_2}\ldots b_r^{k_r}v \right)
 \mod\left( \wp\Omega_F^n + \de \Omega_F^{n-1}\right).
\end{align*}
Using the induction hypothesis on these terms finishes the proof.
\hfill$\square$
\end{prf}

\begin{lem}\label{36}
For  $v\in\Omega_F^{n-1}$, $b_1,\ldots ,b_r\in F$ and $t,k_1,\ldots ,k_r\in\NN$ $t\geq 1$, we can find forms 
$\omega,\omega_1,\ldots,\omega_r\in\Omega_F^{n-1}$ and $q_1,\ldots,q_r\in\NN$ with $0\leq q_1,\ldots,q_r < p^t$ and $k_i\equiv q_i \mod p^t$, so that the following equation holds
\begin{align*}
 b_1^{k_1}\ldots b_r^{k_r} (\de v )^{[p^t]}&\equiv  b_1^{q_1}\ldots b_r^{q_r} (\de \omega )^{[p^t]} + 
 \sum_{i=1}^r b_i \de \omega_i\mod (\wp \Omega_F^n + \de \Omega_F^{n-1})
\end{align*}
\end{lem}

\begin{prf} 
We may assume $b_1,\ldots,b_r\in F^*$ and for at least one $j\in\{1,\ldots,r\}$ we have $k_j\geq p^t$, otherwise there is nothing to prove. 
So assume $k_1\geq p^t$. Note once more that a change of basis changes the image of $s_p$ by an exact form. Using this and computing in $H_p^{n+1}(F)$ we get
\begin{align*}
&\ \ \ \ b_1^{k_1}\ldots b_r^{k_r} (\de v )^{[p^t]} \\
&\equiv b_1^{k_1-p^t}b_2^{k_1}\ldots b_r^{k_r} (b_1 \de v )^{[p^t]} \\
&\equiv b_1^{k_1-p^t}b_2^{k_1}\ldots b_r^{k_r} \left(\de \left( b_1 v \right) \pm v \w \de b_1 \right)^{[p^t]} \\
&\equiv b_1^{k_1-p^t}b_2^{k_1}\ldots b_r^{k_r} \left(\de \left( b_1 v \right)^{[p^t]} \pm b_1^{p^t-1}v^{[p^t]}\w \de b_1   + \de x \right) \\
&\equiv b_1^{k_1-p^t}b_2^{k_1}\ldots b_r^{k_r} \de \left( b_1 v \right)^{[p^t]}    \pm \left( b_1^{k_1-1}b_2^{k_1}\ldots b_r^{k_r} v^{[p^t]} \right) \w \de b_1 
+b_1^{k_1-p^t}b_2^{k_1}\ldots b_r^{k_r} \de x \\
&\equiv b_1^{k_1-p^t}b_2^{k_1}\ldots b_r^{k_r} \de \left( b_1 v \right)^{[p^t]}    +   b_1 \de \left(  \pm b_1^{k_1-1}b_2^{k_1}\ldots b_r^{k_r} v^{[p^t]} \right) 
 + b_1^{k_1-p^t}b_2^{k_1}\ldots b_r^{k_r} \de x \\
 &\quad\mod (\wp\Omega_F^n + \de \Omega_F^{n-1})
\end{align*}
for some $x\in \Omega_F^{n-1}$. This computation can be repeated until the exponent of $b_1$ is smaller than $p^t$ and 
the last summand can be rewritten using Lemma \ref{366}. When this is done, we can make the same steps to lower the exponents of $b_2,\ldots b_r$
if necessary, to finish the proof.
\hfill$\square$
\end{prf}

\begin{defi}
Let $b_1,\ldots ,b_r \in F$ and $m_1,\ldots ,m_r \geq 1$ be integers with $m:=\max\{m_1,\ldots,m_r\}$. 
We define the set $\mathcal{K}_F((b_1,m_1),\ldots,(b_r,m_r))$ to be the subgroup of $H_p^{n+1}(F)$ which is 
additively generated by the forms
\begin{enumerate}
\item [(i)]   $\overline{ b_{i}\de z } $,\ \   with  $z\in \Omega_F^{n-1}$ , $i=1,\ldots,r$  
 \item[(ii)]  $\overline{ b_{1}^{k(t,1)}\cdot\ldots\cdot b_{r}^{k(t,r)} \left( \de v \right)^{[p^t]}   }$, \ \  
 with  $v\in\Omega_F^{n-1}\ , \   t=1,\ldots,m-1$   and \\ $0\leq k(t,i)< p^t$  with  $\max\{1,p^{t-m_i+1}\} \big| k(t,i)$ for $i=1,\ldots,r$.
\end{enumerate}
Note that for $n=0$, we have $\mathcal{K}_F((b_1,m_1),\ldots,(b_r,m_r))=\{0\}$, since $\de \Omega_F^{-1}=\{0\}$.
\end{defi}

In the upcoming results, we will show that the group $\mathcal{K}_F((b_1,m_1),\ldots,(b_r,m_r))$ is the same as $H_p^{n+1}(F(\sqrt[p^{m_1}]{b_1},\ldots,\sqrt[p^{m_r}]{b_r})/F)$. 
For this reason, we will need some properties of these groups.

\begin{lem}\label{377}
Let $(b_1,m_1),\ldots,(b_r,m_r)\in F^* \times \NN_{>0}$. The following assertions hold
\begin{enumerate}
\item[(i)]		$\mathcal{K}_F((b_1,m_1),\ldots,(b_r,m_r))=\mathcal{K}_F((b_{\sigma(1)},m_{\sigma(1)}),\ldots,(b_{\sigma(r)},m_{\sigma(r)}))$ for all $ \sigma \in S_r$.
\item[(ii)]		$\mathcal{K}_{F}((b_1,m_1),\ldots,(b_i^p,m_i+1),\ldots ,(b_r,m_r)) = \mathcal{K}_{F}((b_1,m_1),\ldots,(b_r,m_r))$ for all 
$i\in \{1,\ldots,r\}$.
\end{enumerate}
\end{lem}

\begin{prf}
Part (i) easily follows from the definition of $\mathcal{K}_F((b_1,m_1),\ldots,(b_r,m_r))$. 

Because of (i), for (ii) it is enough to prove 
$$\mathcal{K}_1:=\mathcal{K}_{F}((b_1^p,m_1+1),(b_2,m_r),\ldots ,(b_r,m_r)) = \mathcal{K}_F((b_1,m_1),\ldots,(b_r,m_r))=:\mathcal{K}_2.$$
Set $m:=\max\{m_1,\ldots,m_r\}$, $m'=\max\{m_1+1,m_2,\ldots,m_r\}$ and let $v\in \Omega_F^{n-1}$ be an arbitrary form throughout the whole proof.
Let us start with $\mathcal{K}_1 \subset \mathcal{K}_2$ for which we will show, that each generator of $\mathcal{K}_1$ can be found in $\mathcal{K}_2$. Let us start with the 
generators of type (i). For $i=2,\ldots, r$ the forms $\overline{b_i \de v }$ lie in $\mathcal{K}_2$. Further $\overline{\de\left(b^pv\right)}=\overline{ b_1^p \de v  }=\overline{0}\in \mathcal{K}_2$ obviously holds as well.

So now let us start with a generator of type (ii) in $\mathcal{K}_1$ given by 
$$\overline{ \left(b_{1}^p\right)^{k(t,1)}\cdot\ldots\cdot b_{r}^{k(t,r)} \left( \de v \right)^{[p^t]}   }$$
 with $1 \leq t \leq m'-1$  and $ 0\leq k(t,i) \leq p^t-1$ with $\max\{1,p^{t-m_i+1}\} \big| k(t,i)$ for $i=2,\ldots,r$ and $\max\{1,p^{t-m_1}\}\big| k(t,1)$.

\textbf{Case 1:} $t=m$ (i.e. $m=m_1$ and $\max\{m_1+1,m_2,\ldots,m_r\}=m+1$)\\
Since $p^m$ is not allowed as exponent of $\de v$ for generators of $\mathcal{K}_2$, we will lower the exponent as follows. Since $\max\{1,p^{m-m_i+1}\} \big| k(m,i)$ for all 
$i=2,\ldots,r$ and $m_i\leq m$ holds, we get $p \big| k(m,i)$ and we may write $k(m,i)=p\ell(m,i)$ for suitable $\ell(m,i)$. Using Lemma \ref{30} we get
\begin{equation}\label{gen1}
\overline{ \left(b_{1}^p\right)^{k(m,1)}\ldots b_{r}^{k(m,r)} \left( \de v \right)^{[p^m]}   } = 
\overline{ b_{1}^{k(m,1)} b_2^{\ell(m,2)}\ldots b_{r}^{\ell(m,r)} \left( \de v \right)^{[p^{m-1}]}   }.
\end{equation}
Now we can apply Lemma \ref{36} on equation \eqref{gen1} to find forms  $\omega,\omega_i \in \Omega_F^{n-1}$ and $s_1\in \NN$ with 
$0\leq k(m,1)-s_1 p^{m-1}<p^{m-1}$ such that 
\begin{align*}
&\quad \ \overline{ b_{1}^{k(m,1)} b_2^{\ell(m,2)}\ldots b_{r}^{\ell(m,r)} \left( \de v \right)^{[p^{m-1}]}   } \\
&=\overline{  b_{1}^{k(m,1)-s_1p^{m-1}} b_2^{\ell(m,2)}\ldots b_{r}^{\ell(m,r)} \left( \de \omega \right)^{[p^{m-1}]}   } + \sum_{i=1}^r \overline{  b_i \de \omega_i  }.
\end{align*}
Note that for $k(m,1)< p^{m-1}$ we have $s_1=0$, $\omega=v$ and $\omega_i=0$, since we do not need to lower the exponent of $b_1$. We will now show that the first term of the right hand side of the 
new equation is a generator of $\mathcal{K}_2$. For $m=1$ the term can be readily rewritten using Lemma \ref{366}. So let us assume $m>1$. 
We will check that the exponent of $b_1$ satisfies the necessary conditions for generators in $\mathcal{K}_1$. Obviously we have 
$0\leq k(m,1)-s_1p^{m-1} < p^{m-1}$ by construction and since $m_1=m$ we have $\max\{1,p^{m-1-m+1}\}=1\big| k(m,1)-s_1p^{m-1}$. Now we check that for $i=2,\ldots,r$ 
the exponents of the $b_i$ satisfy the necessary conditions as well. Because of $0\leq k(m,i)=p\ell(m,i)<p^m$, we readily get $0\leq \ell(m,i) < p^{m-1}$ and the condition 
$\max\{1,p^{m-1-m_i+1}\}=\max\{1,p^{m-m_i}\} \big| \ell(m,i)$ holds since we have $\max\{1,p^{m-m_i+1}\}\big| p\ell(m,i)$ by assumption.

\textbf{Case 2:} $t<m$ \\
Now we have a generator of $\mathcal{K}_1$ given by 
$$\overline{ \left(b_{1}^p\right)^{k(t,1)}\cdot\ldots\cdot b_{r}^{k(t,r)} \left( \de v \right)^{[p^t]}   }
=\overline{ b_{1}^{pk(t,1)}\cdot\ldots\cdot b_{r}^{k(t,r)} \left( \de v \right)^{[p^t]}   }$$
with $1 \leq t \leq m-1$  and $ 0\leq k(t,i) \leq p^t-1$ with $\max\{1,p^{t-m_i+1}\} \big| k(t,i)$ for $i=2,\ldots,r$ and $\max\{1,p^{t-m_1}\}\big|k(t,1)$. 
Applying Lemma \ref{36} on this generator, we 
find forms  $\omega,\omega_i \in \Omega_F^{n-1}$ and $s_1\in \NN$ with 
$0\leq pk(t,1)-s_1 p^{t-1}<p^{t-1}$ such that  
$$ \overline{ b_{1}^{pk(t,1)}\cdot\ldots\cdot b_{r}^{k(t,r)} \left( \de v \right)^{[p^t]} }
=\overline{  b_{1}^{pk(t,1)-s_1p^{t}} b_2^{k(t,2)}\ldots b_{r}^{k(t,r)} \left( \de \omega \right)^{[p^{t}]}   } + \sum_{i=1}^r \overline{  b_i \de \omega_i  }.$$
Note once more, that for $pk(t,1)<p^t$, we do not need to lower the exponent of $b_1$ and may choose $s_1=0$, $\omega=v$ and $\omega_i=0$. We will now check once again, that the 
first term of the right hand side of this equation is a generator of $\mathcal{K}_2$. Like above, for $t=1$ we use Lemma \ref{366} to rewrite this term only using generators of type (i). So now 
assume $t>1$. First we have $1 \leq t \leq m-1$. For the $i=2,\ldots,r$ the conditions 
$0\leq k(t,i)\leq p^t-1$ and 
$\max\{1,p^{t-m_i+1}\} \big | k(t,i)$ are obviously true by the choice of the $k(t,i)$. Last we have to check that $\max\{1,p^{t-m_1+1}\}\big| pk(t,1)-s_1p^t$ holds. This is clearly true for 
$m_1\geq t$. For $m_1< t$, by the choice of $k(t,1)$ we have $\max\{1,p^{t-m_1}\}\big| k(t,1)$, thus we can find $\ell(t,1)$ with $pk(t,1)=p\cdot p^{t-m_1}\ell(t,1)$. With 
this $p^{t-m_1+1}=\max\{1,p^{t-m_1+1}\} \big| pk(t,1)-s_1p^t= p^{t-m_1+1}(\ell(t,1)-s_1p^{m_1-1})$ obviously holds. So now we have $\mathcal{K}_1 \subset \mathcal{K}_2$.

Let us now prove that  $\mathcal{K}_2 \subset \mathcal{K}_1$ is also true. For that we will check once again that each generator of  $\mathcal{K}_2$ can be found in 
$ \mathcal{K}_1$. Let us start with the generators of type (i). For $i=2,\ldots,r$ $\overline{b_i \de v}\in \mathcal{K}_1$ is clearly true. Since 
$\max\{m_1+1,m_2,\ldots,m_r\}\geq 2$, we have $\overline{  b_1^p (\de v)^{[p]}  } = \overline{ b_1 \de v  }\in \mathcal{K}_1$.

So now choose a generator of type (ii) of $\mathcal{K}_2$ given by
\begin{equation}\label{gen2}
\overline{ b_{1}^{k(t,1)}\ldots b_{r}^{k(t,r)} \left( \de v \right)^{[p^t]}   }
\end{equation} 
with $1 \leq t \leq m-1$  and $ 0\leq k(t,i) \leq p^t-1$ with $\max\{1,p^{t-m_i+1}\} \big| k(t,i)$ for $i=1,\ldots,r$.

\textbf{Case 1:} $p \big|  k(t,1)$ \\
In this case we can find $\ell(t,1)\in \NN$ with $k(t,1)=p\ell(t,1)$. Thus we can rewrite \eqref{gen2} as 
$$\overline{ \left(b_{1}^p\right)^{\ell(t,1)} b_2^{k(t,2)}\ldots b_{r}^{k(t,r)} \left( \de v \right)^{[p^t]}   }   $$
and we will now prove that this is a generator of $\mathcal{K}_1$. First we have $1\leq t \leq m-1 \leq m'-1$. For $i=2,\ldots,r$  we have $0\leq k(t,i)\leq p^t-1$ and
$\max\{1,p^{t-m_i+1}\}\big| k(t,i)$ by the choice of the $k(t,i)$. Fir $i=1$ it suffices to check that $\max\{1,p^{t-(m_1+1)+1}\}\big|\ell(t,1)$ is true, 
since $0\leq \ell(t,1)\leq p^t-1$ clearly holds. This condition is true for 
$m_1\geq t$, and for $m_1<t$, we get $k(t,1)=p^{t-m_1+1}q(t,1)$ for some $q(t,1)\in\NN$, thus $\ell(t,1)=p^{t-m_1}q(t,1)$ which is clearly divisible by $p^{t-m_1}=\max\{1,p^{t-m_1}\}$.

\textbf{Case 2:} $p\not\big|k(t,1)$ \\
Applying Lemma \ref{30} to equation \eqref{gen2}, we get 
\begin{equation}\label{gen3}
\overline{ b_{1}^{k(t,1)}\ldots b_{r}^{k(t,r)} \left( \de v \right)^{[p^t]}   } = \overline{ \left(b_{1}^p\right)^{k(t,1)} b_2^{pk(t,2)}\ldots b_{r}^{pk(t,r)} \left( \de v \right)^{[p^{t+1}]}   }.
\end{equation} 
Will will now check, that the right hand side \eqref{gen3} is a generator of $\mathcal{K}_1$. Note that since $p\not\big|k(t,1)$, we have $t<m_1$. 
If $m_1$ is the maximum of $\{m_1,\ldots,m_r\}$, we have $m'=m+1$, thus  
$1\leq t+1 \leq m'-1$ holds. If $m_1<m$ (i.e. $\max\{m_1+1,m_2,\ldots,m_r\}=m'=m$), we get 
$1\leq t+1 \leq m_1\leq m-1$. So in both cases the exponent $p^{t+1}$ of $\de v$ meets the necessary conditions. By the choice of $k(t,i)$ for $i=2,\ldots,r$, 
we have $\max\{1,p^{t-m_i+1}\}\big|k(t,i)$, thus $\max\{1,p^{(t+1)-m_i+1}\}\big| pk(t,i)$ and $0\leq pk(t,i)\leq p^{t+1}-1$ is true. Finally for the exponent of $b_1^p$, the conditions 
$\max\{1,p^{(t+1)-(m_1+1)+1}\}=\max\{1,p^{t-m_1+1}\}\big| k(t,1)$ and $0\leq k(t,1)\leq p^{t+1}-1$ also hold by the choice of $k(t,1)$. This finishes the proof of (ii).

\hfill$\square$

\end{prf}

Before checking that the group $\mathcal{K}_F((b_1,m_1),\ldots,(b_r,m_r))$ is actually the kernel of the extension  $F(\sqrt[p^{m_1}]{b_1},\ldots,\sqrt[p^{m_r}]{b_r})/F$, we need to check 
this claim for the cases discussed in \ref{3.3} and \ref{34}. This will be done in the following Lemma.

\begin{lem}\label{378}
Let $(b_1,m_1),\ldots,(b_r,m_r)\in F^* \times \NN_{>0}$ with $m:=\max\{m_1,\ldots,m_r\}$.
\begin{enumerate}
\item[(i)]	For $[F(\sqrt[p^{m_1}]{b_1},\ldots,\sqrt[p^{m_r}]{b_r}):F]=p$, there exist  $j\in \{1,\ldots ,r\}$ and $c_j\in F\setminus F^p$ such that 
 $b_j = c_j^{p^{m_j-1}}$ and 
$\mathcal{K}_F((b_1,m_1),\ldots,(b_r,m_r))= \overline{ c_j \de \Omega_F^{n-1}   }$.
\item[(ii)]	For $m=1$ we have 
$\mathcal{K}_F((b_1,m_1),\ldots,(b_r,m_r)) = \sum_{i=1}^r  \overline{ b_i \de \Omega_F^{n-1} }$.
\end{enumerate}
\end{lem}

\begin{prf}
(i) Set $E:=F(\sqrt[p^{m_1}]{b_1},\ldots,\sqrt[p^{m_r}]{b_r})$. Because of $[E:F]\neq 1$, there must be a $j\in\{1,\ldots,r\}$ with $b_j\not\in F^{p^{m_j}}$. 
If there would be a $k\in\{1,\ldots,r\}$ with $b_k\not\in F^{p^{m_k-\ell}}$ with $\ell\geq 2$, then we would have $[F(\sqrt[p^{m_k}]{b_k}):F]\geq p^2$ which leads to 
$[E:F]\geq p^2$. Hence we have $b_i\in F^{p^{m_i-1}}$ for all $i=1,\ldots,r$ and there is a $j\in\{1,\ldots,r\}$ with 
$b_j\in F^{p^{m_j-1}}\setminus F^{p^{m_j}}$. Assume there is another $j\neq k\in\{1,\ldots r\}$ with $b_k\in F^{p^{m_k-1}}\setminus F^{p^{m_k}}$. Then write $b_j=c_j^{p^{m_j-1}}$ and 
$b_k=c_k^{p^{m_k-1}}$ with $c_j,c_k \in F \setminus F^p$. If $c_j,c_k$ would be $p$-independent, then $E$ would contain the subfield 
$F(\sqrt[p^{m_j}]{b_j},\sqrt[p^{m_k}]{b_k})=F(\sqrt[p]{c_j}, \sqrt[p]{c_k})$, 
which has degree $p^2$ over $F$, which again leads to $[E:F]\geq p^2$. 

Using Lemma \ref{377}, without loss of generality we may assume $j=1$ and write 
$b_1=c_1^{p^{m_1-1}}$ with $c_1\in F\setminus F^p$. For all $i\neq 1$, one 
of the following cases hold
\begin{itemize}
\item [(a)] $b_i\in F^{p^{m_i}}$, i.e. $b_i=c_i^{p^{m_i}}$ for $c_i\in F$.
\item [(b)] $b_i\in F^{p^{m_i-1}}\setminus F^{p^{m_i}}$ with $b_i=c_i^{p^{m_i-1}}$ and $c_i\in F^p(c_1)$.
\end{itemize}
By using Lemma \ref{377} once more, we may rearrange the $b_i$ so that we can find $s$,  $2\leq s\leq r$, 
such that $b_2,\ldots,b_s$ meet condition (b) and $b_{s+1},\ldots,b_r$ meet condition (a). All that is now 
left to prove is
$$\mathcal{K}:=\mathcal{K}_F((b_1,m_1),\ldots,(b_r,m_r))= \overline{ c_1 \de \Omega_F^{n-1}   }.$$
The inclusion $\overline{ c_1 \de \Omega_F^{n-1}   } \subset \mathcal{K}$ is trivial, since for $v\in\Omega_F^{n-1}$ we have $\overline{  b_1 (\de v)^{[p^{m_1-1}]} }= \overline{  c_1 \de v }\in \mathcal{K}$ 
by Lemma \ref{30}. To prove the reversed inclusion, we will once again show that each generator of $\mathcal{K}$ is contained in $\overline{c_1\de \Omega_F^{n-1}}$ and we will start with generators 
of type (i). For the rest of the proof fix $v\in\Omega_F^{n-1}$. We have 
$$\overline{ b_1 \de v }= \overline{ c_1^{p^{m_1-1}} \de v }= \begin{cases} \overline{ c_1 \de v}  \ &, \ m_1=1 \\ \overline{0} \ &, \ m_1>1  \end{cases}$$
For $i=2,\ldots,s$ write $c_i=\sum_{k=0}^{p-1} x_{ik}^p c_1^k$ for some $x_{ik} \in F$. For $m_i>1$, we get $\overline{ b_i \de v } =0$ and for $m_i=1$ we get
$$\overline{ b_i \de v }= \sum_{k=1}^{p-1}  \overline{c_1^k \de \left( x_{ik}^p v \right)}
 \overset{\text{Lemma } \ref{366}}{=} \sum_{k=1}^{p-1} \overline{c_1 \de \left( k c_1^{k-1} x_{ik}^p v \right)} = \overline{c_1 \de \left( \sum_{k=1}^{p-1} k c_1^{k-1} x_{ik}^p v \right) }.$$
Lastly, if $i=s+1,\ldots,r$ we get $\overline{ b_i \de v } = \overline{ 0 }$ because of $b_i\in F^{p^{m_i}}\subset F^p$.

Now choose a generator of type (ii) given by 
$$\overline{ b_{1}^{k(t,1)}\cdot\ldots\cdot b_{r}^{k(t,r)} \left( \de v \right)^{[p^t]}   }$$
with $1 \leq t \leq m-1$, $0\leq k(t,i)\leq p^t-1$ and $\max \{1,p^{t-m_i+1}\}\big| k(t,i)$ for $i=1,\ldots,r$.
If $m_i=1$ for some $i\in \{1,\ldots,r\}$, we would get $k(t,i)=0$, so without loss of generality $m_i\geq 2$ for $i=1,\ldots,r$. Computing in $H_p^{n+1}(F)$, we get 
\begin{align*}
&\ \quad b_{1}^{k(t,1)}\cdot\ldots\cdot b_{r}^{k(t,r)} \left( \de v \right)^{[p^t]} \\
&\equiv c_1^{p^{m_1-1}k(t,1)}c_2^{p^{m_2-1}k(t,2)}\ldots c_s^{p^{m_s-1}k(t,s)}c_{s+1}^{p^{m_{s+1}}k(t,s+1)}\ldots c_{r}^{p^{m_{r}}k(t,r)} \left( \de v \right)^{[p^t]} \\
&\overset{(1)}{\equiv} c_1^{p^{m_1-1}k(t,1)}c_2^{p^{m_2-1}k(t,2)}\ldots c_s^{p^{m_s-1}k(t,s)} \left( \de \left( c_{s+1}^{p^{m_{s+1}-t}k(t,s+1)}\ldots c_{r}^{p^{m_{r}-t}k(t,r)} v \right) \right)^{[p^t]}\\
&\!\!\!\!\!\overset{\text{Lem } \ref{30}}{\equiv} c_1^{p^{m_1-1-t}k(t,1)}c_2^{p^{m_2-1-t}k(t,2)}\ldots c_s^{p^{m_s-1-t}k(t,s)}  \de \left( c_{s+1}^{p^{m_{s+1}-t}k(t,s+1)}\ldots c_{r}^{p^{m_{r}-t}k(t,r)} v \right)\\
&\mod (\wp \Omega_F^n + \de\Omega_F^{n-1}),
\end{align*}
where step (1) holds because $p^{m_i}k(t,i)$ is divisible by $p^{t+1}$ due to the definition of $k(t,i)$. 
Set $\omega := c_{s+1}^{p^{m_{s+1}-t}k(t,s+1)}\ldots c_{r}^{p^{m_{r}-t}k(t,r)} v$. Since $c_1,\ldots c_s\in F^p(c_1)$, we can write 
$$  c_1^{p^{m_1-1-t}k(t,1)}c_2^{p^{m_2-1-t}k(t,2)}\ldots c_s^{p^{m_s-1-t}k(t,s)} = \sum_{k=0}^{p-1} y_k^p c_1^k$$
for some $y_k \in F$. Using this, we finally see 
\begin{align*}
&\quad\  c_1^{p^{m_1-1-t}k(t,1)}c_2^{p^{m_2-1-t}k(t,2)}\ldots c_s^{p^{m_s-1-t}k(t,s)}  \de \left( c_{s+1}^{p^{m_{s+1}-t}k(t,s+1)}\ldots c_{r}^{p^{m_{r}-t}k(t,r)} v \right) \\
 &\equiv \left( \sum_{k=1}^{p-1} y_k^p c_1^k \right) \de \omega 
\overset{\text{Lemma }\ref{366}}{\equiv} c_1 \de \left( \sum_{k=0}^{p-1} k c_1^{k-1} y_k^p \omega    \right)
\mod (\wp \Omega_F^n + \de\Omega_F^{n-1}),
\end{align*}
which finishes the proof of (i). The assertions (ii) is an easy consequence of the definition of $\mathcal{K}_F((b_1,m_1),\ldots,(b_r,m_r))$, since no generators of type (ii) appear.

\hfill$\square$

\end{prf}

Using the notions from the proof above, we just showed that for $[E:F]=p$, we have $F(\sqrt[p^{m_1}]{b_1},\ldots,\sqrt[p^{m_r}]{b_r})=F(\sqrt[p]{c_j})$, 
thus $\mathcal{K}_F((b_1,m_1),\ldots,(b_r,m_r))=H_p^{n+1}(E/F)$. And if $\max\{m_1,\ldots,m_r\}=1$, obviously $\mathcal{K}_F((b_1,m_1),\ldots,(b_r,m_r))=H_p^{n+1}(E/F)$. With these two 
cases covered, we can now prove the following result.

\begin{theo}\label{1}
Let $E=F(\sqrt[p^{m_1}]{b_1},\ldots,\sqrt[p^{m_r}]{b_r})/F$ be a purely inseparable extension with $b_1,\ldots,b_r \in F^*$ and integers $m_1,\ldots,m_r \geq 1$. Then we have
$$H_p^{n+1}(E/F)=\mathcal{K}_F((b_1,m_1),\ldots,(b_r,m_r)).$$
\end{theo}

\begin{prf}
Set $m:=\max\{m_1,\ldots,m_r\}$ and $\mathcal{K}:=\mathcal{K}_F((b_1,m_1),\ldots,(b_r,m_r))$. 
We will start by showing that each generator of $\mathcal{K}$ is actually an element of $H_p^{n+1}(E/F)$. 
This is clear for forms of the type (i) since $b_{j}\in E^p$ for all $j=1,\ldots,r$. Let us now take a look at forms of the type (ii).   
For this we choose $t\in \{1,\ldots, m-1\}\, , \, v\in \Omega_F^{n-1}$ and for $j=1,\ldots,r$ we take $k(t,j)$ with $0\leq k(t,j)\leq p^t-1$ and $\max\{1,p^{t-m_j+1}\}\big|k(t,j)$. Set
$\beta_j:= \sqrt[p^{m_{j}}]{b_{j}}$. Over $E$, we have 
\begin{align*}
b_{1}^{k(t,1)}\ldots b_{r}^{k(t,r)}(\de v)^{[p^t]} &\equiv \beta_1^{\left(p^{m_{1}}\right)k(t,1)}\ldots \beta_r^{\left(p^{m_{r}}\right)k(t,r)}(\de v)^{[p^t]}  \\
&\overset{(1)}{\equiv}  \left( \beta_1^{\left(p^{m_{1}-t}\right)k(t,1)}\ldots \beta_r^{\left(p^{m_{r}-t}\right)k(t,r)} \de v \right)^{[p^t]} \\
&\!\!\!\!\!\!\!\! \overset{\text{Lemma }\ref{30}}{\equiv} \beta_1^{\left(p^{m_{1}-t}\right)k(t,1)}\ldots \beta_r^{\left(p^{m_{r}-t}\right)k(t,r)} \de v \\
&\overset{(2)}{\equiv}  \de \left(  \beta_1^{\left(p^{m_{1}-t}\right)k(t,1)}\ldots \beta_r^{\left(p^{m_{r}-t}\right)k(t,r)} v\right) \\
&\equiv 0 \mod (\wp\Omega_E^n + \de \Omega_E^{n-1}). 
\end{align*}
Step (1) holds, because for all $j=1,\ldots,r$, the exponent $p^{m_{j}}k(t,j)$ is divisible by $p^t$. This is clear for $m_{j}\geq t+1$ and also for $m_{j}\leq t$ since, 
by definition of $k(t,j)$, we find $\ell(t,j)$ with $k(t,j)=\ell(t,j)p^{t-m_{j}+1}$, hence $p^{m_{j}}k(t,j)=p^{t+1}\ell(t,j)$. Finally step (2) holds, because 
all the exponents $(p^{m_{1}-t})k(t,1),\ldots ,(p^{m_{r}-t})k(t,r)$ are divisible by $p$ by the just given argument.

Let us now prove the reverse inclusion. First assume there is a $j\in\{1,\ldots,r\}$ with $b_j\in F^p$. Then we can find $c_j\in F$ with $b_j=c_j^p$ and we get 
$$E=F(\sqrt[p^{m_1}]{b_1},\ldots,\sqrt[p^{m_j-1}]{c_j},\ldots,\sqrt[p^{m_r}]{b_r}).$$
So changing $(b_j,m_j)\in F^*\times \NN_{>0}$ to $(c_j,m_j-1)$ has no influence on the field extension $E$, hence it does not change the kernel $H_p^{n+1}(E/F)$. By Lemma \ref{377}, 
the group $\mathcal{K}_F((b_1,m_1),\ldots,(b_r,m_r))$ is also fixed under the change $(b_j,m_j)\mapsto (c_j,m_j-1)$. Using this, without loss of generality we may assume 
$b_1,\ldots,b_r\in F\setminus F^p$ and we further assume $m=m_1\geq \ldots \geq m_r$.

Say $[E:F]=p^h$ for some $h\in \NN$. We will prove the theorem by 
induction on $h$.

The case $h=1$ is covered by Lemma $\ref{378}$, so assume $h>1$. For $m=m_1=1$, we have $E=F(\sqrt[p]{b_1},\ldots,\sqrt[p]{b_r})$ and this case was also dealt with in \ref{378}.

Now say $m=m_1\geq 2$. To simplify the following computations, let us fix some notations. Set $B:=b_2\cdot \ldots \cdot b_r$ and for all $\nu=(\nu_2,\ldots,\nu_r)\in \NN^{r-1}$ set 
$B^{\nu}=b_2^{\nu_2}\cdot \ldots \cdot b_r^{\nu_r}$. On the set $\NN^{r-1}$, we will use the usual scalar multiplication and addition, which are both defined component-wise. 
Define $M:=F \left( \sqrt[p]{b_1} \right)\subset E$. With this we get 
$$E=M\left(\sqrt[p^{m_1}]{b_1},\ldots,\sqrt[p^{m_r}]{b_r} \right).$$
Note that some of the $b_j$ might be $p^{\text{th}}$ powers over $M$. Also note that $[E:M]=p^{h-1}$. For $t=1,\ldots,m-1$ set 
\begin{align*}
Q_t:=\{ (k_2,\ldots,k_r)\in \NN^{r-1} \mid\  &0\leq k_j \leq p^t-1 \text{ and } \max\{1,p^{t-m_j+1}\}\big| k_j \\& \text{ for } j=2,\ldots,r\}
\end{align*}
The sets $Q_t$ contain all the possible exponents of the elements $b_2,\ldots,b_r$ for generators of type (ii) as defined in $\mathcal{K}$ for a fixed $t$.

Now choose $\overline{\omega} \in H_p^{n+1}(E/F)$. Thus we have $\overline{\omega_M} \in H_p^{n+1}(E/M)$ and because of $[E:M]=p^{h-1}$, we have
\begin{align}\label{001}
\omega_M = \wp(u) + \de v + \sum_{j=1}^{r} b_j \de x_j + \sum_{t=1}^{m-1}\sum_{\varepsilon \in Q_t}\sum_{k=0}^{p^t-1} b_1^k B^{\varepsilon} (\de z_{t\varepsilon k})^{[p^t]}
\end{align}
with $u\in\Omega_M^n$ and $v,x_j,z_{t\varepsilon k} \in \Omega_M^{n-1}$ for all suitable $j,t,\varepsilon,k$.


Since $b_1\in F\setminus F^p$, we may choose a $p$-basis $\mathcal{B}_F$ of $F$ containing $b_1$ and obtain a $p$-basis of $M$ 
by $\mathcal{B}_M= \left( \mathcal{B}_F\setminus\{b_1\}\right) \cup \{\sqrt[p]{b_1}\}$. 
Thus, using the decomposition of $\Omega_M^n$ described in Remark \ref{04}(ii), for every $i=0,\ldots,p-1$ we can find unique forms 
$u_i \in \Omega_F^n$ , $v_i, x_{ji} , z_{t \varepsilon k i}\in \Omega_F^{n-1}$ and $u'\in M \de\left( \sqrt[p]{b_1}\right) \w \Omega_F^{n-1}$ , 
$v',x_j', z_{t\varepsilon k}'\in M \de\left( \sqrt[p]{b_1}\right) \w \Omega_F^{n-2}$ with
\begin{align*}
u&=\sum_{i=0}^{p-1} \left( \sqrt[p]{b_1}\right)^i u_i  + u' \ \ ,\  \ 
v=\sum_{i=0}^{p-1} \left( \sqrt[p]{b_1}\right)^i v_i  + v'\, , \\
x_j&=\sum_{i=0}^{p-1} \left( \sqrt[p]{b_1}\right)^i x_{ji}  + x_j'  \ \ , \ \ 
z_{t\varepsilon k}=\sum_{i=0}^{p-1} \left( \sqrt[p]{b_1}\right)^i z_{t\varepsilon k i}  + z_{t\varepsilon k}',
\end{align*}
Applying the additive maps $\wp , \de$ and $(s_p)^t$, we obtain
\begin{align*}
\wp(u) &= \sum_{i=0}^{p-1} b_1^i u_i^{[p]} - \left( \sqrt[p]{b_1}\right)^i u_i  + u''      \ \ ,\  \ 
\de v = \sum_{i=0}^{p-1} \left( \sqrt[p]{b_1}\right)^i \de v_i + v'' \\
\de x_j &= \sum_{i=0}^{p-1} \left( \sqrt[p]{b_1}\right)^i \de x_{ji} + x_j''     \ \ ,\  \ 
(\de z_{t\varepsilon k})^{[p^t]}= \sum_{i=0}^{p-1} b_1^{ip^{t-1}} (\de z_{t\varepsilon k i})^{[p^t]}  + z_{t\varepsilon k}''
\end{align*}
with some suitable $u'' \in M \de\left( \sqrt[p]{b_1}\right) \w \Omega_F^{n-1}$ and $v'',x_j'', z_{t\varepsilon k}''\in M \de\left( \sqrt[p]{b_1}\right) \w \Omega_F^{n-2}$.
Inserting this in Equation \eqref{001}, we get
\begin{align*}
\omega_M = &\sum_{i=0}^{p-1}\left( b_1^i u_i^{[p]} - \left( \sqrt[p]{b_1}\right)^i u_i  + u'' \right)
+ \sum_{i=0}^{p-1}\left( \left( \sqrt[p]{b_1}\right)^i \de v_i + v'' \right)\\
&+ \sum_{j=1}^{r}  \sum_{i=0}^{p-1} \left(  b_j\left( \sqrt[p]{b_1}\right)^i \de x_{ji} + b_j x_j'' \right) \\
&+ \sum_{t=1}^{m-1}\sum_{\varepsilon \in Q_t}\sum_{k=0}^{p^t-1} \sum_{i=0}^{p-1} \left(  b_1^{k+ip^{t-1}}B^{\varepsilon} (\de z_{t\varepsilon k i})^{[p^t]}  +  b_1^k B^{\varepsilon}z_{t\varepsilon k}'' \right) 
\end{align*}

By Remark \ref{04}, we have 
$$\Omega_M^n=\bigoplus_{i=0}^{p-1} \left(\sqrt[p]{b_1}\right)^i \left( \Omega_F^n\right)_M \oplus 
\bigoplus_{i=0}^{p-1} \left(\sqrt[p]{b_1}\right)^i \de \left(\sqrt[p]{b_1} \right) \w \left( \Omega_F^{n-1}\right)_M.$$
Using this decomposition on the equation above, we get $2p$ different equation, but we are only interested in those not containing the slot $\de \left(\sqrt[p]{b_1} \right)$. Thus we get
\begin{align}\label{0012}
\omega_M &= \wp(u_0) +\de v_0+ \sum_{i=1}^{p-1} b_1^i u_i^{[p]} 
+ \sum_{j=1}^{r}     b_j \de x_{j0}  \\
&\ \ + \sum_{t=1}^{m-1}\sum_{\varepsilon \in Q_t}\sum_{k=0}^{p^t-1} \sum_{i=0}^{p-1} \left(  b_1^{k+ip^{t-1}}B^{\varepsilon} (\de z_{t\varepsilon k i})^{[p^t]}  \right) \nonumber
\end{align}
and for $i=1,\ldots,p-1$
\begin{align*}
0= - u_i + \de v_i + \sum_{j=1}^r b_j \de x_{ji}.
\end{align*}
So for $i=1,\ldots,p-1$ we have 
\begin{align*}
u_i^{[p]}=  (\de v_i)^{[p]} + \sum_{j=1}^r b_j^p (\de x_{ji})^{[p]}.
\end{align*}
Inserting these in Equation \eqref{0012} and reducing it $\!\! \mod (\wp\Omega_F^n + \de \Omega_F^{n-1})$ we get 
\begin{align*}
 \omega_M    &\equiv \sum_{i=1}^{p-1}  b_1^i(\de v_i)^{[p]}   + \sum_{j=1}^{r}     b_j \de x_{j0}   + 
\sum_{i=1}^{p-1}\sum_{j=1}^r  b_1^i b_j^p (\de x_{ji})^{[p]}    \\
&\ \ + \sum_{t=1}^{m-1}\sum_{\varepsilon \in Q_t}\sum_{k=0}^{p^t-1} \sum_{i=0}^{p-1}    b_1^{k+ip^{t-1}}B^{\varepsilon} (\de z_{t\varepsilon k i})^{[p^t]} 
\ \ \ \mod (\wp\Omega_F^n + \de \Omega_F^{n-1}) . \nonumber
\end{align*}
Even though only containing forms defined over $F$, this equation is still defined over the field $M$. By reducing it to the field $F$, we have to add some 
element of $H_p^{n+1}(M/F)=\overline{  b_1\de \Omega_F^{n-1}  }$. Thus we can 
find a form $y\in\Omega_F^{n-1}$ with 
\begin{align*}
 \omega   &\equiv \sum_{i=1}^{p-1}  b_1^i(\de v_i)^{[p]}   + b_1 \de y + \sum_{j=1}^{r}     b_j \de x_{j0}  
+ \sum_{i=1}^{p-1}\sum_{j=1}^r  b_1^i b_j^p (\de x_{ji})^{[p]}     \\
&\ \ + \sum_{t=1}^{m-1}\sum_{\varepsilon \in Q_t}\sum_{k=0}^{p^t-1} \sum_{i=0}^{p-1}    b_1^{k+ip^{t-1}}B^{\varepsilon} (\de z_{t\varepsilon k i})^{[p^t]} 
\ \ \ \mod (\wp\Omega_F^n + \de \Omega_F^{n-1})  \nonumber
\end{align*}
Using Lemma \ref{36} on the last two sums to lower the exponents whenever needed (i.e. whenever an exponent is larger than $p$ resp. $p^t$), 
we see that $\overline{\omega}$ can be expressed as a sum of generators of $\mathcal{K}$.

\hfill$\square$
\end{prf}

\begin{rem}\label{05} Choose $E=F(\sqrt[p^{m_1}]{b_1},\ldots,\sqrt[p^{m_r}]{b_r})$ like in Theorem \ref{1}.

(i) Note that this generating system of $H_p^{n+1}(E/F)$ is independent of a possible $p$-dependence of the elements $b_1,\ldots,b_r$. But the number of generators can be reduced in 
the case of a non modular extension. For example assume that $b_{i_1},\ldots,b_{i_q}$ is a maximal $p$-independent subset of $b_1,\ldots, b_r$ (i.e. $[F^p(b_1,\ldots,b_r):F^p]=p^q$ and 
$F^p(b_1,\ldots,b_r)=F^p(b_{i_1},\ldots,b_{i_q})$). The forms of type (i) can then be 
replaced by the forms 
\begin{align*}
(i')\ \ \overline{ b_{i_j} \de z } \text{ with }  z\in \Omega_F^{n-1}, \ j=1,\ldots,q
\end{align*}

(ii) Now assume $b_1,\ldots, b_r\in F\setminus F^p$ (which can be assumed with every purely inseparable extension). The exponent $\exp(E/F)$ is then given by $\max \{m_1,\ldots,m_r\}$, which shows that the exponent 
of $E/F$ has direct influence on the generators given in Theorem \ref{1}. Further choose $s\in \{1,\ldots,e\}$ and set $s_i=\min\{s,m_i\}$ for $i=1\ldots,r$. It is easy to see that 
$F(\sqrt[p^{s_1}]{b_1},\ldots,\sqrt[p^{s_r}]{b_r})=\{x\in E \mid x^{p^k}\in F$ for a $k\leq s\}$. With this 
we obviously have
$$\mathcal{K}_F((b_1,s_1),\ldots,(b_r,s_r)) = H_p^{n+1}(F(\sqrt[p^{s_1}]{b_1},\ldots,\sqrt[p^{s_r}]{b_r})/F)$$
and
$$\mathcal{K}_F((b_1,s_1),\ldots,(b_r,s_r)) \subset \mathcal{K}_F((b_1,m_1),\ldots,(b_r,m_r)).$$

(iii) One should also note that not all possible combinations of the $t$ and $k(t,j)$ in Theorem \ref{1} produce different forms. 
For example assume $t\geq 1$ and all the $k(t,j)$ are divisible by $p$, let us say $k(t,j)=p\ell(t,j)$. 
Then we readily get for each $v\in\Omega_F^{n-1}$
\begin{align*}
\overline{  b_{1}^{p\ell(t,1)} \ldots b_{r}^{p\ell(t,r)} (\de v)^{[p^t]}  } 
&=  \overline{\left(  b_{1}^{\ell(t,1)} \ldots b_{r}^{\ell(t,r)}(\de v)^{[p^{t-1}]} \right)^{[p]}  } \\
&\!\!\!\!\!\!\!\!\! \overset{\text{Lemma }\ref{30}}{=}  \overline{b_{1}^{\ell(t,1)} \ldots b_{r}^{\ell(t,r)}(\de v)^{[p^{t-1}]}  }.
\end{align*}

(iv) In the case $p=2$ Theorem \ref{1} was obtained independently by Aravire, Laghribi and O'Ryan using somewhat different methods.
\end{rem}

Remark \ref{05}(iii) actually showed the following

\begin{coro}\label{35}
Let $b\in F\setminus F^p $ and $E=F\left(\sqrt[p^e]{b}\right)$ be a simple purely inseparable field extension of $F$ of exponent $e$. Then
\begin{align*}
H_p^{n+1}(E/F)= \sum_{t=0}^{e-1} \sum_{\substack{1\leq j \leq p^t \\ j,p \textrm{ coprime } }} \overline{  b^j \left( \de\Omega_F^{n-1} \right) ^{[p^t]} }.
\end{align*}
\end{coro}

Our result on kernels $H_p^{n+1}(E/F)$ for finite purely inseparable extensions can be readily extended to possibly infinite purely inseparable extension. For that recall that every purely inseparable extension can be written as $E=F\left( \beta_i \mid i\in I \right)$ with a suitable indexset 
$I$ and elements $\beta_i \in E$. Further we define the exponent of $x\in E$ by $\exp(x/F):=\min\{k\in \NN \mid x^{p^k} \in F   \}$. With this we can easily generalize Theorem 
\ref{1} to the following.

\begin{theo}\label{0}
Let $E/F$ be a purely inseparable extension with $E=F(\beta_i \mid i\in I)$. Set $e_i:=\exp(\beta_i/F)$ and $b_i :=\beta_i^{p^{e_i}}$ for $i\in I$. For a finite subset $J\subset I$, we 
set $e(J):= \exp(F(\beta_j \mid j\in J)/F)=\max\{e_j \mid j\in J\}$. 
The kernel $H_p^{n+1}(E/F)$ is additively generated by the forms 
\begin{enumerate}
\item [(i)] $\overline{ b_{i}\de z } $ \ \  with  $z\in \Omega_F^{n-1}$ , $i\in I$   
 \item[(ii)]$\overline{ \left( \prod_{j\in J} b_j^{k(t,j)} \right) \left( \de v \right)^{[p^t]}   }$
 with  $v\in\Omega_F^{n-1}$, a finite subset $J \subset I$,  $ t=1,\ldots,e(J)-1$  and   
 $0 \leq k(t,j) \leq p^t-1$ with $\max\{1,p^{t-e_j+1}\}\big| k(t,j)$ for $j\in J$.
\end{enumerate}
\end{theo}

\begin{prf}
Choose $\overline{\omega}\in H_p^{n+1}(E/F)$. Since the relation $\overline{\omega}_E=0$ only uses a finite number of elements of $E$, it is clear that there exists already a finite indexset 
$J\subset I$, with $\overline{\omega}\in H_p^{n+1}(F(\beta_j \mid \, j\in J)/F)$. Applying Theorem \ref{1} on the finite extension $F(\beta_j \mid \, j\in J)/F$ finishes the proof.

\hfill$\square$
\end{prf}

Using the notations from Theorem \ref{0}, note that for a finite $J \subset I$ with $e(J)\leq 1$, there are no generators of type (ii).

We want to finish this chapter by finding a different representation for the generators of the kernel $H_p^{n+1}(E/F)$, so they can be used more easily in the later sections.

For that, let us call a generator of  $H_p^{n+1}(E/F)$ of type (ii) of degree $t$, when it is given by  $\overline{ b_{1}^{k(t,1)}\ldots b_{r}^{k(t,r)} ( \de v )^{[p^t]}   }$ 
for suitable $t$ and $k(t,j), j=1,\ldots,r$ and cannot be written with a smaller exponent of $\de v$ (i.e. not all the exponents $k(t,1),\ldots,k(t,r)$ are divisible by $p$). 
By using this notation and combining Theorem \ref{1} with Remark \ref{05}(ii), we get the following useful result.

\begin{coro}\label{37}
Let $E=F(\sqrt[p^{m_1}]{b_1},\ldots,\sqrt[p^{m_r}]{b_r})/F$ be a purely inseparable field extension 
of exponent $e$ over $F$ with $b_1,\ldots,b_r\in F\setminus F^p$ and integers $m_1,\ldots,m_r\geq 1$.
For every $1\leq t \leq e-1$ and $k_1,\ldots,k_r \in \NN$ with $\max\{ 1,p^{t-m_{j}+1} \} \Big| k_j$ for $j=1,\ldots,r$, the form 
$\overline{  b_{1}^{k_1}\cdot\ldots\cdot b_{r}^{k_r} \left( \de v \right)^{[p^t]}   }$
lies in $H_p^{n+1}(E/F)$ and can be written as a sum of generators of $H_p^{n+1}(E/F)$ of degree $\ell$ with $\ell\leq t$.
\end{coro}

With this corollary we can rewrite the generators of $H_p^{n+1}(E/F)$ with $E=F(\sqrt[p^{m_1}]{b_1},\ldots,\sqrt[p^{m_r}]{b_r})/F$ for some $b_1,\ldots,b_r\in F\setminus F^p$. 
Let us start with a 
generator of type (ii), which is given by $\overline{ b_{1}^{k(t,1)}\cdot\ldots\cdot b_{r}^{k(t,r)} \left( \de v \right)^{[p^t]}   }$ 
for suitable $t$ and $k(t,j)$. Since the group $\de \Omega_F^{n-1}$ is additively generated 
by forms of the type $\de a_1 \w \ldots \w \de a_n$ with $a_1\ldots,a_n\in F$ and both $s_p$ and $\de\,$ are additive maps, it suffice to 
use a form $\de a_1 \w \ldots \w \de a_n$ instead of $\de v$.
Note once again that the map $s_p$ is basis-dependent, but a change of basis only changes the image by an exact form. Knowing this, in $\Omega_F^n$ we have
\begin{align*}
 &\ \ \ \   \left( \de a_1 \w \ldots \w \de a_n \right) ^{[p^t]} \\
 &= \left( a_1\ldots a_n \frac{\de a_1}{a_1} \w \ldots \w \frac{\de a_n}{a_n} \right) ^{[p^t]} \\
 &= \left( (a_1\ldots a_n)^p \frac{\de a_1}{a_1} \w \ldots \w \frac{\de a_n}{a_n}+  \de x_{t-1} \right) ^{[p^{t-1}]} \\
  &=   \left( (a_1\ldots a_n)^p \frac{\de a_1}{a_1} \w \ldots \w \frac{\de a_n}{a_n} \right) ^{[p^{t-1}]} 
  +   \left(\de x_{t-1} \right) ^{[p^{t-1}]} \\
&\qquad\vdots \\
&=  \left( (a_1\ldots a_n)^{p^t} \frac{\de a_1}{a_1} \w \ldots \w \frac{\de a_n}{a_n} \right) 
+\sum_{h=0}^{t-1}  \left(\de x_{h} \right) ^{[p^{h}]} \\
&\overset{(3)} {=}  \left( (a_1\ldots a_n)^{p^t} 
\frac{\de\left( a_1\ldots a_n \right)}{(a_1\ldots a_n)}\w\frac{\de a_2}{a_2} \w \ldots \w \frac{\de a_n}{a_n} \right) 
+\sum_{h=0}^{t-1}   \left(\de x_{h} \right) ^{[p^{h}]} 
\end{align*}
for some suitable $x_0\ldots,x_{t-1}\in\Omega_F^{n-1}$. Step (3) in the above computation was made by using the equality 
$\frac{\de x}{x}\w\frac{\de y}{y}=\frac{\de \, ( xy ) }{xy}\w\frac{\de y}{y}$ repeatedly. Multiplying the equation above by $b_{1}^{k(t,1)}\cdot\ldots\cdot b_{r}^{k(t,r)}$, we get 
\begin{align*}
&\qquad b_{1}^{k(t,1)}\cdot\ldots\cdot b_{r}^{k(t,r)} \left( \de a_1 \w \ldots \w \de a_n \right) \\
&\equiv  b_{1}^{k(t,1)}\cdot\ldots\cdot b_{r}^{k(t,r)}\left( (a_1\ldots a_n)^{p^t} 
\frac{\de\left( a_1\ldots a_n \right)}{(a_1\ldots a_n)}\w\frac{\de a_2}{a_2} \w \ldots \w \frac{\de a_n}{a_n} \right) \\
&\qquad+\sum_{h=0}^{t-1} b_{1}^{k(t,1)}\cdot\ldots\cdot b_{r}^{k(t,r)} \left(\de x_{h} \right) ^{[p^{h}]}  
\qquad \mod (\wp \Omega_F^n + \de \Omega_F^{n-1}).
\end{align*}
 Because of Corollary \ref{37}, we can rewrite each form 
$ \overline { b_{1}^{k(t,1)}\cdot\ldots\cdot b_{r}^{k(t,r)}\left(\de x_{h} \right) ^{[p^{h}]}  }$ 
for $h=0,\ldots,t-1$ as sum of generators of degree equal to $h$ and below. So using induction on $t$ and replacing 
$a_1\ldots a_n=:s$, we just showed that the generators of type (ii) can be replaced by the forms 
\begin{align*}
\overline{  b_{1}^{k(t,1)}\cdot\ldots\cdot b_{r}^{k(t,r)}  s^{p^t} \frac{\de s}{s} \w  \frac{\de a_2}{a_2} \w \ldots \w \frac{\de a_n}{a_n}   }
\end{align*}
with $s,a_2,\ldots,a_n\in F$ , 
$t=1,\ldots,e-1$ and for $j=1,\ldots,r$ , $1\leq k(t,j)\leq p^t-1$ with $\max \{1,p^{t-m_{j}+1} \} \Big| k(t,j)$. 
A similar computation can be done with the generators of type (i), which leads to the following result.

\begin{coro}\label{101}
Let $E=F(\sqrt[p^{m_1}]{b_1},\ldots,\sqrt[p^{m_r}]{b_r})/F$ be a purely inseparable extension of exponent $e$ with $b_1,\ldots,b_r \in F\setminus F^p$ and integers $m_1,\ldots,m_r\geq 1$. 
For $n \geq 1$, the kernel $H_p^{n+1}(E/F)$ is additively generated by the forms 
\begin{enumerate}
\item [(i)] $\overline{  b_{j} s \frac{\de s}{s} \w \frac{\de a_2}{a_2} \w \ldots \frac{ \de a_n}{a_n}  }$\ \  with $j= 1,\ldots,r$ and $s,a_2,\ldots,a_n \in F$
 \item[(ii)] $ \overline{ b_{1}^{k(t,1)}\cdot\ldots\cdot b_{r}^{k(t,r)}  s^{p^t} \frac{\de s}{s} \w  \frac{\de a_2}{a_2}\ldots \w \frac{\de a_n}{a_n}  }$\ 
with $s,a_2,\ldots,a_n\in F$, 
$t=1,\ldots,e-1$ and for $j=1,\ldots ,r$ , $ 0 \leq k(t,j) \leq p^t-1$ with $\max \{1,p^{t-m_{j}+1} \} \Big| k(t,j)$.
\end{enumerate}
For $n=0$, we have $H_p^1(E/F)=\{\overline{0}\}$.
\end{coro}

Since this computation was independent of $n$, we also get

\begin{coro}\label{38}
For the purely inseparable extension $E=F(\sqrt[p^{m_1}]{b_1},\ldots,\sqrt[p^{m_r}]{b_r})/F$ with $b_1,\ldots,b_r\in F\setminus F^p$ and integers $m_1,\ldots,m_r\geq 1$, 
we have
\begin{align*}
H_p^{n+1}(E/F)= H_p^2(E/F) \w \overline{ \nu_{n-1}(F) }
\end{align*}
\end{coro}

\end{section}

\begin{section}{The kernels $\Omega^n(E/F)$ and $\nu_n(K/F)$ for a modular purely inseparable field extensions}\label{sec4}

Since we just have studied the behavior of the group $H_p^{n+1}(F)$ under purely inseparable field extension, we will now continue by studying  $\Omega_F^n$
under some special cases of these kind of extensions. For that recall that we have a decomposition given by 
$$ \Sigma_n=\Sigma_{n,>s}\mathbin{\dot{\cup}}\Sigma_{n,\leq s}$$
described in Remark \ref{04}(i). We will omit the index $n$ in the sequel, since it will be fixed during all our computations.

\begin{theo}\label{41}
Let $E=F\left( \sqrt[p^{m_1}]{b_1},\ldots,\sqrt[p^{m_r}]{b_r} \right)/F$ be a purely inseparable field extension of $F$ with $p$-independent $b_1,\ldots,b_r \in F$ 
and integers $m_1,\ldots,m_r\geq 1$. Then
$$\Omega^n(E/F) = \sum_{i=1}^r \de b_i \w \Omega_F^{n-1}. $$
In particular, we have
\begin{align*}
\nu_n(E/F)&= \{ \omega_1 \w \frac{\de b_1}{b_1} + \ldots + \omega_r \w \frac{\de b_r}{b_r} \in \nu_n(F)\ \mid\   \omega_1,\ldots,\omega_r \in \Omega_F^{n-1}      \}.
\end{align*}
\end{theo}

\begin{prf}
Because of $b_j \in E^p$ for all $j=1,\ldots,r$, we have $(\de b_j)_E=0$ which readily shows the first inclusion. So let us now prove the reverse inclusion. Since $b_1,\ldots,b_r$ are $p$-independent, we may 
choose a $p$-basis $\mathcal{B}_F=\{c_i \mid i\in I\}$ of $F$ 
with $c_{i_1}=b_1,\ldots ,c_{i_r} =b_r$ where $i_1,\ldots i_r \in I$ are the first $r$ elements in $I$ with $i_1<\ldots <i_r$. With this, we get a $p$-basis 
of $E$ by $\mathcal{B}_E=\{c_i \mid i\in I\setminus\{i_1,\ldots,i_r \} \} \cup \{ \sqrt[p^{m_1}]{b_1},\ldots,\sqrt[p^{m_r}]{b_r} \}$.

Choose $\omega\in \Omega_F^n$. Thus we can find $x_{\varepsilon},y_{\tau}\in F$ with
\begin{equation}\label{eq40}
\omega = \sum_{\varepsilon\in \Sigma_{>r}} x_{\varepsilon} \frac{\de c_{\varepsilon}}{c_{\varepsilon}} + 
\sum_{\tau \in \Sigma_{\leq r}} y_{\tau} \frac{\de c_{\tau}}{c_{\tau}}
\end{equation}
Since each form $\frac{\de c_{\tau}} {c_{\tau}}$ with $\tau \in \Sigma_{\leq r}$ contains at least one of the forms $\de b_1 , \ldots, \de b_s$, 
which all become zero in $E$, we readily get
$$\omega_E =0= \sum_{\varepsilon\in \Sigma_{>r}} x_{\varepsilon} \left(\frac{\de c_{\varepsilon}}{c_{\varepsilon}}\right)_E .$$
Using the fact that the forms $\left(\frac{\de c_{\varepsilon}}{c_{\varepsilon}}\right)_E$ are part of a basis of $\Omega_E^n$, we get $x_{\varepsilon}=0$ for all 
$\varepsilon\in \Sigma_{>r}$. Inserting this in Equation \eqref{eq40}, we see
\begin{align*}
\omega =\sum_{\tau \in \Sigma_{\leq r}} y_{\tau} \frac{\de c_{\tau}}{c_{\tau}} \in \sum_{j=1}^r \de b_j \w \Omega_F^n
\end{align*}

\hfill$\square$
\end{prf}

The fact that this kernel is independent of the $m_1,\ldots,m_r$ is not surprising at all, since in \cite[Theorem 5.2]{Hof5} Hoffmann showed a similar result for the 
isotropy behavior of bilinear forms, which is in some way similar to the behavior of differential forms. (see the next section)

In the last section we saw that for $E=F\left( \sqrt[p^{m_1}]{b_1},\ldots,\sqrt[p^{m_r}]{b_r} \right)$ a possible $p$-dependence of the elements $b_1,\ldots,b_r$ has no major influence 
on the generators of the kernel $H_p^{n+1}(E/F)$. For $\Omega^n(E/F)$ however this is not the case, which can be seen by the following example.

\begin{bsp}
Choose $p=2$, set $F:=\FF_2(X,Y,Z)$ with some variables $X,Y,Z$ and define the field extension $E/F$ of exponent 2 as $E=F\left( \sqrt[4]{Z},\sqrt[4]{X^2Z+Y^2} \right)$. 
Note that this is a non modular extension. It is commonly known that eight is the smallest possible degree of a purely inseparable field extension in 
characteristic 2, such that there exists a non modular extension. Set $n=2$ and assume 
$$\Omega^2(E/F)=\de Z \w \Omega_F^{1} + \de \left( X^2Z+Y^2 \right)\w \Omega_F^{1} = \de Z \w \Omega_F^{1}.$$
Since $X\sqrt{Z}+Y$ and $\sqrt{Z}$ are squares in $E$, we have 
$$0 = \de \left(  X\sqrt{Z}+Y \right)= \sqrt{Z} \de X + \de Y, $$
which leads to $(\de X \w \de Y)_E=\sqrt{Z} \de X \w \de X=0$, thus $\de X \w \de Y \in \Omega^2(E/F)$. But since $X,Y,Z$ are $2$-independent over $F$, 
we obviously have $\de X \w \de Y \not\in \de Z \w \Omega_F^{1}$.
\end{bsp}

\end{section}

\begin{section}{Quadratic and bilinear forms in characteristic 2}\label{sec6}

In this last section, we will translate the kernels we computed to quadratic and bilinear forms in characteristic two. To do so, we will start by recapping some basic facts about these forms.
We refer to \cite{Elm1} for any undefined terminology or any other facts that we do not mention explicitly. 

In this section, $F$ will always denote a field of \cha 2, if not stated otherwise. By $W(F)$ resp. $W_q(F)$, we denote the Witt ring of non 
degenerated symmetric bilinear forms over $F$
 resp. the Witt group of nonsingular quadratic forms over $F$. The group $W_q(F)$ is endowed with a $W(F)$-modul structure 
by the following multiplication: For a symmetric bilinear form $(\bb, V)\in W(F)$ and a nonsingular quadratic form $(\p, W)\in W_q(F)$, we set $\bb\otimes\p$ to be the 
quadratic form on $V\otimes W$, defined on generators by $\bb\otimes \p(v\otimes w):=\bb(v,v)\p(w)$ for $v\in V, w\in W$ with polar form $\bb_{\bb\otimes\p}=\bb\otimes\bp$. 

Let $I(F)$ denote the fundamental ideal of even dimensional bilinear forms of $W(F)$. 
For an integer $n\geq 1$, set $I^n(F):=\left( I(F)\right)^n$ and let $I^n W_q(F)$ denote the group $I^n(F)\otimes W_q(F)$. 
It is well known that $I^n(F)$ resp. $I^n W_q(F)$ 
is generated by $n$-fold bilinear Pfister forms $\lrp a_1,\ldots , a_n \rrr_{\bb}:= \langle 1,-a_1 \rangle_{\bb} \otimes  \ldots \otimes \langle 1,-a_n \rangle_{\bb}$ resp. by the 
$(n+1)$-fold quadratic Pfister forms $\lrp a_1,\ldots , a_n , b \rrp := \lrp a_1,\ldots,a_n \rrr_{\bb} \otimes [1,b]$, where $\langle c_1,\ldots,c_n \rangle_{\bb}$ 
denotes the symmetric diagonal bilinear form 
$\sum_{i=1}^n c_i X_i Y_i$ and $[e,f]$ denotes the quadratic form $eX^2+XY+fY^2$. By $\HH_{\bb}$ resp. $\mathbb{M}_a$ we denote the bilinear hyperbolic plane resp. the 
bilinear metabolic plane which represents $aF^2$, and by $\HH$ we denote the quadratic hyperbolic plane.

A given field extension $L/F$ induces natural homomorphisms $\iota_L : W(F) \to W(L)$ resp. $\iota_L : W_q(F) \to W_q(L)$, whose kernel $W(L/F)$ is called the bilinear Witt kernel of the extension $L/F$ resp. whose kernel $W_q(L/F)$ is called the quadratic Witt kernel of the extension $L/F$. 

To transport the results from the previous section to bilinear and quadratic forms, we will use the following famous result due to Kazuya Kato in 1982.

\begin{theo}\label{61}{{\rm(\cite{Kat1}, p. 494)}}
Let $F$ be a field of \cha 2. For every $n\in\NN$ there exists isomorphisms
\begin{align*}
e_n : I^n(F) / I^{n+1}(F) &\to \nu_n(F) \\
\lrp a_1,\ldots , a_n \rrr_{\bb} \!\!\!\! \mod I^{n+1}(F) &\mapsto   \frac{\de a_1}{a_1} \w \ldots \w \frac{\de a_n}{a_n}   
\end{align*}
and\begin{align*}
f_n : I^n W_q(F) / I^{n+1} W_q(F) &\to H_2^{n+1}(F) \\
\lrp a_1,\ldots , a_n , b \rrp \!\!\!\! \mod I^{n+1} W_q(F) &\mapsto \overline{  b \frac{\de a_1}{a_1} \w \ldots \w \frac{\de a_n}{a_n}  }
\end{align*}
\end{theo}

Note that $f_0$ can be chosen as the Arf-Invariant, since we have $W_q(F)/I^1 W_q(F) \cong F/\wp(F) = H_2^1(F)$. To simplify notations, 
for every field extension $L/F$ we define the following
\begin{itemize}
\item $\overline{ I^n}(F) := I^n(F)/I^{n+1}(F)$ \quad and \quad   $\overline{ I^n}(L/F) := \ker \left( \overline{ I^n}(F) \to  \overline{ I^n}(L) \right)$
\item $\overline{ I_q^n}(F) := I^nW_q(F)/I^{n+1}W_q(F)$ \quad and \quad   $\overline{ I_q^n}(L/F) := \ker \left( \overline{ I_q^n}(F) \to  \overline{ I_q^n}(L) \right)$
\end{itemize}

We will start by transferring the kernel $H_2^{n+1}(E/F)$ for a finite purely inseparable field extension. Since we already found a suitable representation of the generators of this kernel 
ans since we know how $f_n$ acts on generators, Theorem \ref{61} readily leads to the following result.

\begin{coro}\label{39}
Let $E=F(\sqrt[2^{m_1}]{b_1},\ldots,\sqrt[2^{m_r}]{b_r})/F$ be a purely inseparable field extension 
of exponent $e$ with $b_1,\ldots,b_r\in F\setminus F^2$ and integers $m_1,\ldots,m_r\geq 1$. For $n\geq 1$ the kernel $ \overline{I_q^n}(E/F)$ is 
additively generated by the forms 
\begin{enumerate}
\item[(i)] $\overline{ \lrp a_1,\ldots , a_{n-1} \rrr \otimes  \lrp s, s  b_{j}  \rrp } $ with $ a_ 1,\ldots , a_{n-1}, s\in F$ and $j \in \{1,\ldots,r\} $
\item[(ii)] $\overline{ \lrp a_1,\ldots , a_{n-1} \rrr \otimes  \lrp s, s^{2^t} b_{1}^{k(t,1)}\ldots b_{r}^{k(t,r)} \rrp } $  with  $ a_ 1,\ldots , a_{n-1}, s\in F$, 
  $1 \leq t \leq e-1$ ,  and for  $j=1,\ldots, r \ \ 0\leq k(t,j)\leq 2^t-1$  and  $\max\{1,2^{t-m_{j}+1}\} \Big | k(t,j) $    
\end{enumerate}
For $n=0$, this kernel is trivial.
\end{coro}

\begin{prf}
This proof is done using Corollary \ref{101} together with the commutativity of the diagram 
\begin{align*}
\begin{matrix}
\overline{ I_q^n }(F) & \longrightarrow & \overline{ I_q^n }(E) \\
\downarrow & & \downarrow \\
H_2^{n+1}(F)& \longrightarrow & H_2^{n+1}(E)
\end{matrix}
\end{align*}\hfill$\square$
\end{prf}

We are now able to compute the Witt-kernel for the discussed field extension by following the same steps as Laghribi did in \cite[p.4]{Lag1}. 
For the sake of completeness, we will follow these ideas once again.

\begin{theo}\label{2}
Let $E=F(\sqrt[2^{m_1}]{b_1},\ldots,\sqrt[2^{m_r}]{b_r})/F$ be a purely inseparable field extension 
of exponent $e$ with $b_1,\ldots,b_r\in F\setminus F^2$ and integers $m_1,\ldots,m_r\geq 1$. 
The quadratic Witt-kernel $W_q(E/F)$ is then generated (as a $W(F)$-module) by the forms
\begin{enumerate}
\item[(i)] $\lrp s, s b_{j} ]]$ with $s\in F$ and $j \in\{1,\ldots,r\}$,
\item[(ii)]$ \lrp s, s^{2^t} b_{1}^{k(t,1)}\ldots b_{r}^{k(t,r)} \rrp $
with $s\in F$ ,  $1\leq t\leq e-1$  and for $j=1,\ldots,r\ \\ 0\leq k(t,j)\leq 2^t-1$ and $\max \{1,2^{t-m_{j}+1}\}\ \big|\  k(t,j)$.
\end{enumerate}
In particular, for all $n\geq 1$ we have
\begin{align*}
\ker\left( I^nW_q(F) \to I^nW_q(E) \right) = I^{n-1}(F) \otimes \mathcal{Q}=I^{n-1}(F)\otimes W_q(E/F)
\end{align*} 
where $\mathcal{Q}$ is additively generated by the forms of type (i) and (ii).
\end{theo}

\begin{prf}

Let $G$ be the set generators described in the theorem and $[G]_F$ the set of forms, generated by forms of $G$ over the field $F$ in $W_q(F)$ as a $W(F)$ module. 
First we have to show that each form in $G$ is actually an element of $W_q(E/F)$. Let us start with forms of the type (ii). For this we set $B:=b_{1}\ldots b_{r}$ and 
$\Theta:=\sqrt[2^{m_{1}}]{b_{1}}\ldots \sqrt[2^{m_{r}}]{b_{r}}$. For $\alpha=(\ell_1,\ldots,\ell_r)\in \NN^r$, we define $B^{\alpha}:=b_{1}^{\ell_1}\ldots b_{r}^{\ell_r}$ 
resp. $\Theta^{\alpha}:=\left( \sqrt[2^{m_{1}}]{b_{1}}\right)^{\ell_1}\ldots\left(\sqrt[2^{m_{r}}]{b_{r}}\right)^{\ell_r}$. For $s\in F\ , \ t\in\{1,\ldots,e-1\}$ and 
$\delta := (k(t,1),\ldots,k(t,r))\in \NN^r$ as described in the theorem, we have to prove that the form 
$$\lrp s, s^{2^t}B^{\delta}\rrp$$
becomes hyperbolic over the field $E$. Due to the assumption, the $r$-tupel $\gamma := (2^{m_1} k(t,1),\ldots, 2^{m_r} k(t,r))$ can be written as $\gamma = 2^{t+1} \varepsilon$
with $\varepsilon\in \NN^r$. Using this
together with the fact that the forms $[1,x]$ and $[1,x^{2^{u}}]$ are isometric for all $u\in \NN$ by comparing the Arf-Invariant (see \cite{Elm1}), 
over the field $E$ we get
\begin{align*}
\lrp s, s^{2^t}B^{\delta}\rrp &=\lrp s, s^{2^t}\Theta^{2^{t+1}\varepsilon}\rrp 
=[1,(s \Theta^{2\varepsilon})^{2^{t}}] \perp s [1,(s \Theta^{2\varepsilon})^{2^{t}}]\\
&=[1,s \Theta^{2\varepsilon}] \perp s[1,s \Theta^{2\varepsilon}]
=[\Theta^{2\varepsilon},s] \perp [s, \Theta^{2\varepsilon}] 
= \HH \perp \HH.
\end{align*}
Hence $\lrp s, s^{2^t}B^{\delta}\rrp_E$ is hyperbolic. To see that forms of type (i) also become hyperbolic over $E$ can be proven in a similar way.

Now  let $\p$ be an arbitrary element of $W_q(E/F)$. Without loss of generality, we may choose $\p$ to be 
anisotropic. Using analogous arguments as in the case of the determination of the kernels $H_p^{n+1}(E/F)$, we may assume that $F$ has a finite $2$-basis 
$\mathcal{B}$ with $|\mathcal{B}|=N<\infty$ (see Remark \ref{00}).

Using the basis of $\Omega_F^n$ introduced on page \pageref{basis}, it is easy to see that we have 
$\Omega_F^{i}=0$ for $i>N$, thus $H_2^{i+1}(F)=\{\overline{0}\}$ for $i>N$. Using the isomorphism of Kato, we then have $I^iW_q(F) = I^{i+1}W_q(F)$ for all $i>N$ and by the 
Arason-Pfister-Hauptsatz \cite[4.2]{Bae3}, it follows that 
$I^iW_q(F)=0$ for  $i>N$. Now choose a $k\geq 1$ with $\p \in I^kW_q(F)$. This is possible, because of $\p_E=0$, the form $\p_E$ has a trivial Arf-Invariant over $E$. From
Lemma \ref{355} we see that $\p$ has a trivial Arf-Invariant over $F$ as well, hence $\p \in I W_q(F)$.
So we may assume $\p \in I^k W_q(E/F)$, which leads to $\overline{ \p } \in  \overline{I_q^k }(E/F)$.
Using Corollary \ref{39} and setting $G=\{g_i \mid i\in I\}$, we can find a finite subset $I_0\subset I$ and bilinear forms $\bb^0_{i_0}\in I^{k-1}(F)$ for $i_0\in I_0$ with 
$$\p + \sum_{i_0 \in I_0} \bb^0_{i_0}\otimes g_{i_0} \in I^{k+1}W_q(F).$$
The form $\p + \sum_{i_0 \in I_0} \bb^0_{i_0}\otimes g_{i_0}$ is obviously hyperbolic over the field $E$, since both $\p$ and the $g_{i_0}$ are hyperbolic over $E$. 
If it also hyperbolic over the field $F$, we are finished. If not, 
we can repeat this step to find another finite subset $I_1 \subset I$ and forms $\bb_{i_1}^1\in I^{k}(F)$ for $i_1 \in I_1$ with 
$$\p + \sum_{i_0 \in I_0} \bb^0_{i_0}\otimes g_{i_0} + \sum_{i_1 \in I_1} \bb^1_{i_1}\otimes g_{i_1}  \in  I^{k+2}W_q(F).$$
We can repeat this step until we find finite subsets $I_f \subset I$ and forms $\bb_{i_f}^f \in I^{k-1+f}(F)$ with $i_f \in I_f$ for $f =0,\ldots, N-k$ with 
$$\p + \sum_{f=0}^{N-k} \sum_{i_f \in I_f} \bb^f_{i_f}\otimes g_{i_f}  \in  I^{N+1}W_q(F)=\{0\}.$$
Hence we have 
$$\p = \sum_{f=0}^{N-k} \sum_{i_f \in I_f} \bb^f_{i_f}\otimes g_{i_f}\in [G]_F.$$ 
The second statement can easily be seen by slightly changing the above proof.

\hfill$\square$
\end{prf}

Following the same ideas as in the proof of Theorem \ref{0}, we easily generalize Theorem \ref{2} to possibly infinite purely inseparable field extensions. 

\begin{theo}\label{4}
Let $E/F$ be a purely inseparable extension with $E=F(\beta_i \mid i\in I)$. Set $e_i:=\exp(\beta_i/F)$ and $b_i :=\beta_i^{2^{e_i}}$ for $i\in I$. 
For a finite subset $J\subset I$, we 
set $e(J):= \exp(F(\beta_j \mid j\in J)/F)=\max\{e_j \mid j\in J\}$. 
The quadratic Witt kernel $W_q(E/F)$ is then generated (as a $W(F)$-module) by the forms
\begin{enumerate}
\item[(i)] $\lrp s, s b_{i} ]]$ with $s\in F$ and $i \in I$,
\item[(ii)]$ \lrp s, s^{2^t} \prod_{j\in J} b_j^{k(t,j)} \rrp $
with $s\in F$ ,a finite subset $J \subset I$,   $1\leq t\leq e(J)-1$ and for $j\in J, \ \ 0\leq k(t,j)\leq 2^t-1$ with $\max \{1,2^{t-e_j+1}\}\ \Big|\  k(t,j)$.
\end{enumerate}
In particular, for all $n\geq 1$ we have
$$\ker \left( I^nW_q(F) \to I^nW_q(E) \right) = I^{n-1}(F) \otimes \mathcal{Q'}=I^{n-1}(F)\otimes W_q(E/F)$$
where $\mathcal{Q'}$ is additively generated by the forms of type (i) and (ii).
\end{theo}

\begin{rem}
(i) For possibly infinite purely inseparable extensions of exponent 1, the quadratic Witt-kernel described in Theorem \ref{4} was computed by Hoffmann in \cite{Hof4} only using the theory 
of bilinear and quadratic forms.

(ii) In \cite{Hof3} Hoffmann and Sobiech described the generators of  the quadratic Witt kernel of
the field extension $F(\sqrt[4]{b})/F$ to be 
\begin{enumerate}
\item[(a)] $ \lrp b,c ]]$ with $c\in F^*$, 
\item[(b)] $ \lrp e^3,\frac{b}{e^2} ]]$ with $e\in F^*$.
\end{enumerate}
Setting $f:=e^{-1}$, one readily sees that the generators of type (b) can be rewritten as 
$ \lrp f^{-3},b f^2 ]] \cong \lrp f ,b f^2 ]] $ with $f\in F$, which matches with the generators form Theorem \ref{2}. A generator of type (i) is by Theorem \ref{2} given by 
$\lrp s,bs \rrp$ for some $s\in F$. Since we have 
\begin{align*}
\lrp s,sb \rrp =[1,sb] \perp s[1,sb]= [1,sb] \perp [s,b] =[1,sb]\perp b[1,sb] = \lrp b, sb \rrp,
\end{align*}
we see that these generators are the same as the ones of type (a) by setting $c:=sb$.
\end{rem}

We now want use the kernel $\nu_n(E/F)$ for a modular purely inseparable extension $E/F$ in a similar way to get the kernel for bilinear forms. To do so, we first have to rewrite the generators 
of $\nu_n(E/F)$ in a way that fits the map $e_n$ described in Theorem \ref{61}. This is done by the following result.

\begin{lem}\label{635}
Let $b_1,\ldots,b_r \in F$ be $2$-independent. Then we have
$$\left(\sum_{i=1}^r \de b_i \w \Omega_F^{n-1} \right)\cap \nu_n(F)=\left \langle \frac{\de x}{x}\ \mid x \in F^2(b_1,\ldots,b_r)^* \right \rangle \w \nu_{n-1}(F)$$
\end{lem}

\begin{prf}
The set on the right hand side is obviously additively generated by the logarithmic forms $\frac{\de x}{x}\w \frac{\de y_2}{y_2}\w \ldots \w \frac{\de y_n}{y_n}$ with 
$y_2,\ldots,y_n\in F^*$ and $x\in F^2(b_1,\ldots,b_r)^*$. Now choose one of these generators. Then we can find $\lambda_1,\ldots,\lambda_r\in F$ with 
$\de x = \lambda_1 \de b_1 + \ldots + \lambda_r\de b_r$. With this we readily get 
\begin{align*}
\frac{\de x}{x}\w \frac{\de y_2}{y_2}\w \ldots \w \frac{\de y_n}{y_n} = \sum_{i=1}^r \de b_i \w \left( \lambda_i x^{-1} \frac{\de y_2}{y_2}\w \ldots \w \frac{\de y_n}{y_n} \right).
\end{align*}

Let us now check the reverse inclusion. Choose $\omega\in \left(\sum_{i=1}^r \de b_i \w \Omega_F^{n-1} \right)\cap \nu_n(F)$. Since $b_1,\ldots,b_r$ are $2$-independent, we can extent these elements to a $2$-basis of $F$ given by $\mathcal{B}=\{c_i \mid i\in I\}$. 
We further assume that the elements $b_1,\ldots,b_r$ are the first $r$ elements of $\mathcal{B}$ with ordering $b_1<\ldots < b_r$. Note that by Remark \ref{00}, we may assume this 
$2$-basis to be finite. Using the $F$-basis of $\Omega_F^n$ corresponding to this $2$-basis, we have 
\begin{align}\label{eq66}
\omega= \sum_{\varepsilon \leq \delta} a_{\varepsilon}\frac{\de c_{\varepsilon}}{c_{\varepsilon}}
=a_{\delta}\frac{\de c_{\delta}}{c_{\delta}} +\sum_{\varepsilon < \delta} a_{\varepsilon}\frac{\de c_{\varepsilon}}{c_{\varepsilon}} 
=a_{\delta}\frac{\de c_{\delta}}{c_{\delta}} + \omega_{<\delta},
\end{align}
with $\omega_{<\delta}\in \Omega_{F,<\delta}^n$ and each $\de c_{\varepsilon}$ contains at least one of the forms $\de b_1, \ldots , \de b_r$ by assumption. 
Since $\omega\in \nu_n(F)$, we have $\wp(\omega)\in \de \Omega_F^{n-1}$ and because $\wp$ maintains the filtration of $\omega$, we get 
$$\wp(a_{\delta}\frac{\de c_{\delta}}{c_{\delta}}) \in \de \Omega_F^{n-1} + \Omega_{F,<\delta}^n.  $$
Applying Kato's Lemma \ref{02} (note that $F^{p-1}=F$ for $p=2$), we get 
$$a_{\delta}\frac{\de c_{\delta}}{c_{\delta}} = \frac{\de x_1}{x_1}\w \ldots \w \frac{\de x_n}{x_n} + v_{<\delta}$$
with $v_{<\delta} \in \Omega_{F,<\delta}^n$ and $x_i \in F^*_{\delta(i)}=F^2\left(c_j \mid j\leq \delta(i)\right)^*$. Since $b_1,\ldots,b_r$ are the first $r$ elements of $\mathcal{B}$, 
we have $c_{\delta(1)}\in \{b_1,\ldots,b_r\}$. Thus we get $x_1 \in F^*_{\delta(i)} \subset F^2(b_1,\ldots,b_r)^*$. Inserting this result in Equation \eqref{eq66}, 
we get 
$$\omega = a_{\delta}\frac{\de c_{\delta}}{c_{\delta}} + \omega_{<\delta} = \frac{\de x_1}{x_1}\w \ldots \w \frac{\de x_n}{x_n} + v_{<\delta} +\omega_{<\delta}.$$
Since we already know that $\omega,\frac{\de x_1}{x_1}\w \ldots \w \frac{\de x_n}{x_n}\in \left(\sum_{i=1}^r \de b_i \w \Omega_F^{n-1} \right)\cap \nu_n(F)$ holds, we 
get $v_{<\delta} +\omega_{<\delta}\in \left(\sum_{i=1}^r \de b_i \w \Omega_F^{n-1} \right)\cap \nu_n(F)$ as well. Using induction on 
$ v_{<\delta} +\omega_{<\delta}$ finishes the proof.

\hfill$\square$
\end{prf}

With this lemma, we easily get

\begin{theo}\label{63}
Let $E=F(\sqrt[2^{m_1}]{b_1},\ldots,\sqrt[2^{m_r}]{b_r})/F$ be a purely inseparable field extension with 
$2$-independent $b_1,\ldots,b_r\in F$ and integers $m_1,\ldots,m_r\geq 1$. Then we have 
$$\nu_1(E/F)=\left \langle \frac{\de x}{x}\mid x \in F^2(b_1,\ldots,b_r)^* \right \rangle = \left \{  \frac{\de x}{x}\ \mid x \in F^2(b_1,\ldots,b_r)^* \right\}$$
and 
$$\nu_n(E/F)=\nu_1(E/F)\w \nu_{n-1}(F).$$
\end{theo}

If we now use Kato's isomorphism \ref{61} on $\nu_n(E/F)$ with its new system of generators and follow steps similar to those 
that were made in the proof of Theorem \ref{2}, we get our final result.

\begin{theo}\label{3}
Let $E=F(\sqrt[2^{m_1}]{b_1},\ldots,\sqrt[2^{m_r}]{b_r})/F$ be a purely inseparable field extension with 
$2$-independent $b_1,\ldots,b_r\in F$ and integers $m_1,\ldots,m_r\geq 1$. 
\begin{enumerate}
\item[(i)] $\overline{I^n}(E/F)= I^{n-1}(F) \otimes \overline{I^1}(E/F)$, and generated as an additive subgroup of $\overline{I^1}(F)$ we have
\begin{align*}
\overline{I^1}(E/F) &= \left \langle \overline{  \lrp x \rrr_{\bb}  } \mid  x \in F^2(b_1,\ldots,b_r)^* \right\rangle 
\end{align*}

\item[(ii)] For the bilinear Witt-kernel $W(E/F)$, generated as an ideal of $W(F)$, we have 
\begin{align*}
W(K/F) &= \Big\langle \langle 1,x \rangle_{\bb} \mid x \in F^2(b_1,\ldots,b_r)^* \Big\rangle = \sum_{x \in F^2(b_1,\ldots,b_r)^* } W(F) \otimes \lrp x \rrr_{\bb}
\end{align*}

\item[(iii)] For each $n\in \NN$ we have 
\begin{align*} 
\ker \left( I^n(F) \to I^n(E)  \right) &= I^{n-1}(F) \otimes \left \langle \lrp x \rrr_{\bb} \mid x \in F^2(b_1,\ldots,b_r)^* \right\rangle\\
&=I^{n-1}(F)\otimes W(E/F)
\end{align*}
\end{enumerate}

\end{theo}

We omit the proof, since it is similar to the one from \ref{2}.

\begin{rem}
(i) Theorem \ref{3} can be easily generalized to infinite purely inseparable extension $E/F$, where $E$ is generated over $F$ by roots of $2$-independent elements. This is left to the reader.

(ii) The kernel given in Theorem \ref{3}(ii) was also computed by Hoffmann in \cite{Hof5}, only using the theory of bilinear forms. 

(iii) As a final application we want to note that, using the same notations as in \ref{3}, we have $\nu_1(E/F) \cong \overline{I^1}(E/F) \cong \ker\left( F/F^2 \to E/E^2 \right)$, where 
the second isomorphism is done using the determinant. 
So for every $f\in F\setminus F^2$, we have $f\in E^2$ if and only if $f\in F^2(b_1,\ldots,b_r)$.

\end{rem}

\end{section}

\section*{Acknowledgments and Notes}
The results contained in this paper are part of the PhD-Thesis of the author. I thank D. Hoffmann for supervising this work and giving some very useful hints and corrections in the 
process.

I also thank R. Aravire for some interesting discussions on this topic during his stay at the TU Dortmund in 2014 and for sketching a proof of Theorem \ref{3.3} for the case $p=2$. 
We also thank him for some very helpful remarks on  \cite{Ara1}. R. Aravire's stay was supported by DFG Projekt HO 4784/1-1 {\it Annnihilators and kernels in Kato's cohomology in positive
characteristic and in Witt groups in characteristic 2}.

\bibliographystyle{apalike}
\bibliography{literatur}

\end{document}